\documentclass[reqno]{amsart} \usepackage{amsthm,amssymb}
\usepackage{graphicx}

\def\jb#1{\langle#1\rangle} \def\norm#1{\|#1\|}
\def\normb#1{\bigg\|#1\bigg\|} 
\newcommand{\les}{\lesssim} 
\newcommand{\wt}{\widetilde}

\newcommand{\B}{\mathcal{B}} 
 \newcommand{\F}{\mathcal{F}}
 
\newcommand{\cS}{\mathcal{S}} \newcommand{\cH}{\mathcal{H}}
\newcommand{\cP}{\mathcal{P}}

\newcommand{\T}{\mathbb{T}} \newcommand{\C}{\mathbb{C}}
\newcommand{\N}{\mathbb{N}} \newcommand{\R}{\mathbb{R}}
\newcommand{\Z}{\mathbb{Z}}

\newcommand{\al}{\alpha} \newcommand{\be}{\beta}
\newcommand{\ga}{\gamma} \newcommand{\de}{\delta}
\newcommand{\e}{\varepsilon} \newcommand{\fy}{\varphi}
 
\newcommand{\te}{\theta} \newcommand{\s}{\sigma}
  \newcommand{\x}{\xi}
\newcommand{\y}{\eta}  
 \newcommand{\ft}{{\mathcal{F}}}

\newcommand{\De}{\Delta} \newcommand{\Om}{\Omega}

\newcommand{\p}{\partial} \newcommand{\na}{\nabla}
\newcommand{\re}{\mathop{\mathrm{Re}}}

\newcommand{\supp}{\operatorname{supp}}

\newcommand{\lec}{\lesssim} \newcommand{\gec}{\gtrsim}
 
\newcommand{\IN}[1]{\text{ in }#1}  \newcommand{\etc}{,\ldots,} \newcommand{\I}{\infty}

\newcommand{\ti}{\widetilde} 
 
\newcommand{\LR}[1]{{\langle #1 \rangle}}

\newcommand{\EQ}[1]{\begin{equation}\begin{split}
      #1 \end{split}\end{equation}} 
\newcommand{\EQN}[1]{\begin{equation*}\begin{split}
      #1 \end{split}\end{equation*}} 
 \newcommand{\Del}[1]{}
\newcommand{\BR}[1]{\left\{#1\right\}}
\newcommand{\CAS}[1]{\begin{cases} #1 \end{cases}}

\newcommand{\pt}{&} \newcommand{\pr}{\\ &} \newcommand{\pq}{\quad}
\newcommand{\pn}{}  \newcommand{\prQ}{\\
  &\qquad} 

\numberwithin{equation}{section}

\newtheorem{thm}{Theorem}[section] 
\newtheorem{lem}[thm]{Lemma} \newtheorem{prop}[thm]{Proposition}
\theoremstyle{remark} \newtheorem{rem}[thm]{Remark}

\begin{document}

\title[4D Zakharov system]{Well-posedness and scattering for the
  Zakharov system in four dimensions} \author[I.~Bejenaru, Z.~Guo,
S.~Herr, K.~Nakanishi]{Ioan Bejenaru, Zihua Guo, Sebastian Herr, Kenji
  Nakanishi} \address{Department of Mathematics, University of
  California, San Diego, La Jolla, CA 92093-0112, USA}
\email{ibejenaru@math.ucsd.edu}

\address{School of Mathematical Sciences, Monash University, VIC 3800, Australia \& LMAM, School of Mathematical Sciences, Peking University,
  Beijing 100871, China} \email{zihua.guo@monash.edu}

\address{Fakult\"{a}t f\"{u}r Mathematik, Universit\"{a}t Bielefeld,
  Postfach 10 01 31, 33501 Bielefeld, Germany}
\email{herr@math.uni-bielefeld.de}

\address{Department of Mathematics, Kyoto University, Kyoto 606-8502,
  Japan} \email{n-kenji@math.kyoto-u.ac.jp}

\begin{abstract}
  The Cauchy problem for the Zakharov system in four dimensions is
  considered. Some new well-posedness
  results are obtained.  For small initial data, global well-posedness
  and scattering results are proved, including the case of initial
  data in the energy space.  None of these results is restricted to
  radially symmetric data.
\end{abstract}

\maketitle


\section{Introduction and Main Results}\label{sect:i-main}
Let $\al>0$. The Zakharov system
\begin{equation}\label{eq:Zak}
  \begin{cases}
    i\dot u - \De u = nu,\\
    \ddot n/\al^2 - \De n = -\De|u|^2,
  \end{cases}
\end{equation}
with initial data
\EQ{
 u(0,x)=u_0,\, n(0,x)=n_0,\,\dot n(0,x)=n_1,}
is considered as a simplified mathematical model for Langmuir waves in
a plasma, which couples the envelope $u:\R^{1+d}\to\C$ of the electric
field and the ion density $n:\R^{1+d}\to \R$, neglecting magnetic
effects and the vector field character of the electric field, see
\cite[Chapter V]{SuSu} and \cite{Zak}.

The parameter $\al>0$ is called the ion sound speed.  Formally, as
$\alpha\to \infty$, \eqref{eq:Zak} reduces to the focusing cubic
Schr\"odinger equation
\begin{align}\label{eq:Schr}
  i\dot u - \De u = |u|^2u,
\end{align}
which is energy-critical in dimension $d=4$, see for example
\cite{KM,KV,Dodson} and the references therein concerning recent
developments on global-wellposedness, blow-up and scattering for
\eqref{eq:Schr}. For rigorous results on the subsonic limit (as
$\alpha\to \infty$) of \eqref{eq:Zak} to \eqref{eq:Schr} we refer the
reader to \cite{SW,OT,MN}.

Strong solutions $(u,n)$ of the Zakharov system preserve the mass
\begin{equation}\label{eq:ma}
  \int_{\R^d}|u|^2dx=\int_{\R^d}|u_0|^2dx,
\end{equation}
and the energy, with $D:=\sqrt{-\De}$,
\begin{equation}\label{eq:en}
  E(u,n,\dot n)=\int_{\R^d} |\na u|^2 + \frac{|D^{-1}\dot
    n|^2}{2\al^2}+\frac{|n|^2}{2}-n|u|^2 dx=E(u_0,n_0,n_1).
\end{equation}
In view of \eqref{eq:en}, a natural space for the initial data is the
energy space \EQ{\label{eq:indata-d} (u_0,n_0,n_1)\in H^1(\R^d)\times
  L^2(\R^d)\times \dot H^{-1}(\R^d).}

For initial data in the energy space, the Zakharov system is known to
be globally well-posed if $d=1$ (see \cite{GTV}) and locally well-posed if $d=2,3$
(see \cite{BoCo}). A low regularity local well-posedness theory has been
developed in \cite{GTV} in all dimensions, with further extensions in
\cite{BHHT} if $d=2$, and in \cite{BeHe} if $d=3$, see also the
references therein for previous work. In the case of the torus $\T^d$
well-posedness results were proved in \cite{Takaoka,Kishi}.

In \cite{Merle} blow-up results in finite or infinite time for initial
data of negative energy were proved if $d=3$, and if $d=2$ blow-up in
finite time was derived in \cite{GM1,GM2}. Concerning the final data problem in weighted Sobolev spaces, we refer to \cite{Shimo,GV,OT2}.

Recently, the asymptotic behaviour as $t\to\I$ for the initial data problem was studied in dimension $d=3$: In \cite{GN}, small data
energy scattering in the radial case was obtained by using a normal form technique
and the improved Strichartz estimates for radial functions from \cite{GuoWang}. In \cite{GNW}, a dichotomy between scattering and grow-up was obtained for radial solutions with energy
below the ground state energy. In the non-radial case in dimension
$d=3$, scattering was obtained in \cite{HPS} under the assumption that
the initial data are small enough and have sufficient regularity and
decay. This result was improved recently in \cite{GLNW,Guo}, where
scattering was shown for small initial data belonging to the energy
space with some additional angular regularity.

In the present paper, we continue the analysis of the initial value problem
\eqref{eq:Zak} and focus on the energy-critical dimension $d=4$.  In particular, we will
address the small data global well-posedness and scattering problem in the energy
space, i.e.\ \EQ{\label{eq:indata} (u_0,n_0,n_1)\in H^1(\R^4)\times
  L^2(\R^4)\times \dot H^{-1}(\R^4),}
with no additional symmetry or decay assumption.

We reduce the wave equation to first order equation as usual: Let
\EQ{N:=n - iD^{-1}\dot n/\al,} then $n=\re N=(N+\bar N)/2$ and the
Zakharov system for $(u,N)$ reads as follows:
\begin{equation}\label{eq:Zak1}
  \begin{cases}
    (i\p_t-\De) u =Nu/2+\bar N u/2,\\
    (i\p_t+\al D) N = \al D|u|^2.
  \end{cases}
\end{equation}
The Hamiltonian then becomes
\begin{equation}\label{eq:en-z}
  E(u,n,\dot n)=E_Z(u,N):=\int_{\R^4}|\na u|^2 + \frac{|N|^2}{2}-\re N|u|^2 dx.
\end{equation}
We will restrict ourselves to the system \eqref{eq:Zak1}.
Our first main result is a small data global well-posedness and scattering result.
\begin{thm}\label{thm1-1}
There exists $\e_0=\e_0(\al)>0$ such that for any $(s,l)\in\R^2$
satisfying
 \EQ{ \label{sl range}
 \pt l \ge 0, \pq s<4l+1, \pq (s,l)\not=(2,3),
 \pr \max(\frac{l+1}{2},l-1)\leq s\leq \min(l+2,2l+\frac{11}{8}), }
and for any initial data $(u_0,N_0)\in H^s(\R^4)\times H^l(\R^4)$ satisfying
\EQ{ \label{smallness}
  \norm{(u_0,N_0)}_{H^{1/2}(\R^4)\times L^2(\R^4)}< \e_0,}
there exists a unique global solution $(u,N)\in C(\R;H^s(\R^4)\times H^l(\R^4))$ of \eqref{eq:Zak1} with some space-time integrability. The solution map is continuous in the norms
\EQ{
 H^s \times H^l \ni (u_0,N_0) \mapsto (u,N) \in L^\I(\R;H^s\times H^l).}
Moreover, there exist $(u^\pm, N^\pm)\in H^s(\R^4)\times H^{l}(\R^4)$ such that
\EQ{
 \lim_{t\to \pm \I}(\norm{u(t)-S(t)u^\pm}_{H^s}+\norm{N(t)-W_\al(t)N^\pm}_{H^l})=0,}
where $S(t)=e^{-it\Delta}$ and $W_\al(t)=e^{it\al D}$ are the free propagators.
\end{thm}
In the above statement, we need the space-time integrability to ensure the uniqueness. For example, for any $T>0$,
\EQ{ \label{condof uniq}
 u\in L^2((0,T);B^{1/2}_{4,2}(\R^4))}
is sufficient for uniqueness on $[0,T]$, where $B^{1/2}_{4,2}$ is the inhomogeneous Besov space. See Propositions \ref{small data scatter}, \ref{regup Xs} and \ref{regup higher s} for more detail of the space-time integrability.

Very recently, we learned about an independent work of Kato and Tsugawa \cite{KaTs}.
By a different method, they prove the small data scattering for $l=s-1/2 \ge 0$, using bilinear estimates in $U^p$-$V^p$ spaces for the standard iteration. While their iteration scheme is more direct, our estimates are more elementary and we cover a wider range of $(s,l)$.

Our second result is a large data local well-posedness result for the same range of regularity $(s,l)$ as above, except for the energy space $H^1(\R^4)\times L^2(\R^4)$.
\begin{thm}\label{thm1-2}
Let $(s,l)\in\R^2$ satisfy \eqref{sl range} and $(s,l)\ne (1,0)$. Then, for any $(u_0,N_0)\in H^s(\R^4)\times H^l(\R^4)$, there exists $T=T(u_0,N_0)>0$ and a unique local solution $(u,N)\in C([-T,T];H^s(\R^4)\times H^l(\R^4))$ to \eqref{eq:Zak1} satisfying some space-time integrability {\rm(}\eqref{condof uniq} is enough for the uniqueness{\rm )}.
Both $T>0$ and $(u,N)$ depend continuously on $(u_0,N_0)$.
\end{thm}

\begin{figure}[h]
  \begin{minipage}[c]{.5\textwidth}
    \vspace*{\fill} \centering
    \includegraphics{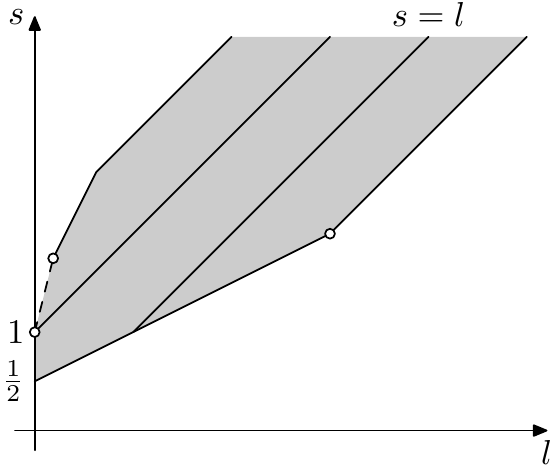}
    \caption{Range of $(s,l)$ obtained in Theorems \ref{thm1-1} and \ref{thm1-2}.}\label{fig}
  \end{minipage}%
  \begin{minipage}[c]{.5\textwidth}
    \vspace*{\fill} \centering
    \includegraphics{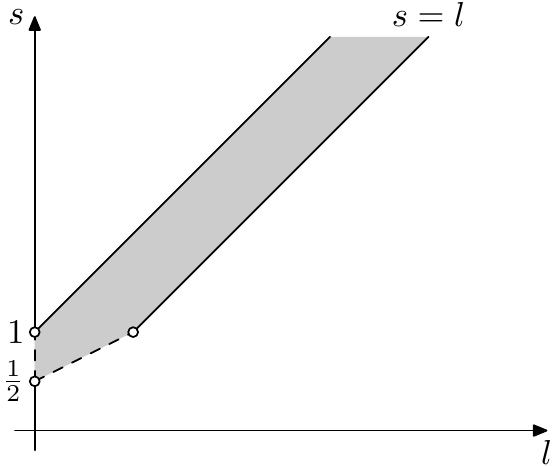}
    \caption{Range of $(s,l)$ obtained in \cite{GTV}.}\label{fig2}
  \end{minipage}
\end{figure}

In dimension $d=4$, Ginibre-Tsutsumi-Velo \cite{GTV} proved local well-posedness in the range $l\leq s \leq l+1$, $l>0, 2s>l+1$, see Figure \ref{fig2}.
Their method is the standard Picard iteration argument in the $X^{s,b}$ spaces.
Theorem \ref{thm1-2} gives further local well-posedness results in some new region indicated in Figure \ref{fig}, while Theorem \ref{thm1-1} covers the same range of exponents as well as the energy space $(s,l)=(1,0)$, which is missing from the large data result Theorem \ref{thm1-2}.

The proofs for Theorems \ref{thm1-1} and \ref{thm1-2} use the normal form
technique and Strichartz estimates as in \cite{GN} and the follow-up papers \cite{GNW,GLNW,Guo}
and related work on the Klein-Gordon-Zakharov system \cite{gnw-kg1,gnw-kg2}.
Our argument is somewhat simpler than \cite{GTV} and it also implies some scattering results.

There is a qualitative difference in our proof between $s<l+1$ and $s>l+1$.
Since the Strichartz norm of $W_\al(t)$ is worse than that of $S(t)$,
we use only the $H^l_x$ norm for $N$, while keeping the full Strichartz norm for $u$, for $s<l+1$.
For $s>l+1$, however, this strategy is prevented by the normal form of $u$, so we need to modify the Strichartz norm for $u$, and to use that of $N$.
Consequently, we can not recover all the Strichartz norms of $S(t)$ for $u$, in spite of the scattering. See Proposition \ref{regup higher s} for the precise statement.
This is consistent with that \cite{GTV} is restricted to $s\le l+1$ and $X^{s,b}$ implies the full range of Strichartz norm.

The energy space $(s,l)=(1,0)$ is at the intersection of $s=l+1$ and
$l=0$, where our multilinear estimates actually break down. More
precisely, we can not close any Strichartz bound for the normal form
of $u$ when $(s,l)=(1,0)$. This is why $(1,0)$ is excluded from
Theorem \ref{thm1-2}. Fortunately enough, with the help of the
conservation law \eqref{eq:en-z} and using the well-posedness in
nearby $(s,l)$, we are still able to show global well-posedness and
scattering in the energy space $(s,l)=(1,0)$ for small data as in
Theorem \ref{thm1-1}. Since the limit NLS \eqref{eq:Schr} is
critical in the energy space $H^1(\R^4)$, it may have blow-up with
bounded $H^1\times L^2$ norm for large data, which suggests that
there may be essential difference between large and small data.

At the other excluded endpoint $(s,l)=(2,3)$, we can prove a strong ill-posedness result, both by instant exit and by non-existence.
\begin{thm} \label{thm-ill}
There exists a radial function $u_0\in H^2(\R^4)$ such that for any $\e>0$, any $N_0\in H^3(\R^4)$, and any $T_0>0$, the system \eqref{eq:Zak1} has no solution $(u,N)\in C([0,T_0];\cS'(\R^4)^2)$ satisfying $(u(0),N(0))=(\e u_0,N_0)$, the equation \eqref{eq:Zak1} in the distribution sense, and
\EQ{ \label{nonexis cond}
 (u,N)\in L^2((0,T_0);H^1(\R^4)\times H^3(\R^4)).}
Moreover, the unique local solution $(u,N)\in C([-T,T];H^2\times H^2)$ given by Theorem \ref{thm1-2} satisfies $N(t)\not\in H^3(\R^4)$ for all $t\in[-T,T]\setminus\{0\}$.
\end{thm}
Note that \eqref{nonexis cond} is weaker than the usual weak solutions, as it does not require $(u(t),N(t))\in H^s\times H^l$ for all $t$ near $0$. The above ill-posedness is due to the mismatch of regularity between $u$ and $N$ in the normal form for $N$.

The rest of paper is organized as follows.
In Section 2, we recall the normal form reduction from \cite{GN}, and then gather multilinear estimates used in the later sections. They easily follow from the Littlewood-Paley decomposition, Coifman-Meyer bilinear estimate, Strichartz and Sobolev inequalities.
Using these estimates and the standard contraction argument, we first prove the small data scattering in $H^s\times H^l$ for $s\le l+1$ in Section 3, and then the local well-posedness for large data in $H^{1/2}\times L^2$ in Section 4. In Section 5, we extend these results to higher regularity by persistence of regularity, except for the energy space $(s,l)=(1,0)$. Theorem \ref{thm1-1} for $(s,l)\not=(1,0)$ follows from Propositions \ref{small data scatter}, \ref{regup Xs} and \ref{regup higher s}.
Similarly, Theorem \ref{thm1-2} follows from Propositions \ref{LWP large}, \ref{regup Xs} and \ref{regup higher s}.
 In Section 6, we prove Theorem \ref{thm1-1} in the energy space $(s,l)=(1,0)$, using the results in $(s,0)$ for $s<1$ and in $(1,l)$ for $l>0$. In Section 7, we prove the ill-posedness Theorem \ref{thm-ill} at $(s,l)=(2,3)$.

\section{Normal form and multilinear estimates}\label{sect:est}
In this section, we set up integral equations and basic estimates for solving the equation. Our analysis is based on the normal form reduction devised in \cite{GN}.
\subsection{Review of the normal form reduction and notation from
  \cite{GN}}\label{subsect:nf-n}
Let $\hat\phi=\F\phi$ denote the Fourier transform of $\phi$. We use
$S(t)$ and $W_\alpha(t)$ to denote the Schr\"odinger and wave
semigroup, respectively:
\[S(t)\phi=\ft^{-1}(e^{it|\xi|^2}\hat{\phi}),\quad
W_\alpha(t)\phi=\ft^{-1}(e^{i\alpha t|\xi|}\hat{\phi}).\] Fix a
radial, smooth, bump function $\eta_0: \R^4 \rightarrow [0, 1]$ with
support in the ball $B_{\frac{8}{5}}(0)$, which is equal to $1$ in the
smaller ball $B_{\frac{4}{5}}(0)$. For $k\in \Z$ let
$\chi_k(\xi)=\eta_0(\xi/2^k)-\eta_0(\xi/2^{k-1})$ and $\chi_{\leq
  k}(\xi)=\eta_0(\xi/2^k)$, and let $P_k, P_{\leq k}$ denote the
corresponding Fourier multipliers.

For two functions $u$, $v$ and a fixed $K \in \N$, $K\geq 5$, we define the paraproduct type operators
\EQ{
  (uv)_{LH}:=&\sum_{k\in\Z}(P_{\leq k-K}u)(P_kv),
  \pq (uv)_{HL}:=(vu)_{LH},\\
  (uv)_{HH}:=&\sum_{\substack{|k_1-k_2|\leq K-1 \\ k_1,k_2\in\Z}}(P_{k_1}u)(P_{k_2}v),}
so that $uv=(uv)_{LH}+(uv)_{HL}+(uv)_{HH}$. We also define
\EQ{
  (uv)_{\al L}:=&\sum_{\substack{|k-\log_2\al|\le 1,\\ k\in\Z}}(P_ku)(P_{\leq k-K}v), \pq (uv)_{L\al}:=(vu)_{\al L},\\
  (uv)_{X L}:=&\sum_{\substack{|k-\log_2\al|> 1,\\ k\in\Z}}(P_ku)(P_{\leq k-K}v),\pq (uv)_{LX}:=(vu)_{X L},}
so that $(uv)_{HL}=(uv)_{\al L}+(uv)_{XL}$.

Moreover, for any $*=HH,LH,HL,\al L$, etc., we denote the symbol (multiplier) of the bilinear operator $(uv)_*$ by $\cP_*$.
We denote finite sums of these bilinear operators in the obvious way, e.g.~$(uv)_{LH+HH}=(uv)_{LH}+(uv)_{HH}$.
With these notations, it was shown in \cite{GN} that \eqref{eq:Zak1} is equivalent --at least for smooth solutions-- to the following integral equation.
Henceforth, for simplicity, we replace the nonlinear term $\re N u/2$ with $N u$ as in \cite{GN}, because the complex conjugation here makes no essential difference for our arguments.
\begin{equation}\label{eq:intu}
  \begin{split}
    u(t)=&S(t)u_0-S(t)\Omega(N,u)(0)+\Omega(N,u)(t)\\
    &-i\int_0^tS(t-s)\Omega(\al D|u|^2,u)(s)ds-i\int_0^tS(t-s) \Omega(N,Nu)(s)ds\\
    &-i\int_0^tS(t-s)(Nu)_{LH+HH+\alpha L}(s)ds,
  \end{split}
\end{equation}
and
\begin{equation}\label{eq:intN}
  \begin{split}
    N(t)=&W_\alpha(t)N_0-W_\alpha(t)D\tilde\Omega(u,u)(0)+D\tilde\Omega(u,u)(t)\\
    &-i\int_0^tW_\alpha(t-s)\al D(u\bar
    u)_{HH+\alpha L+L\al}ds\\
    &-i\int_0^tW_\alpha(t-s)(D\tilde\Omega(Nu,u)+D\tilde\Omega(u,Nu))(s)ds,
  \end{split}
\end{equation}
where $\Omega$, $\tilde \Omega$ are the bilinear Fourier multiplication
operators
\begin{align*}
  \Omega(f,g)=&\F^{-1}\int\cP_{XL}\frac{\hat f(\x-\y)\hat g(\y)}{-|\x|^2+\al|\x-\y|+|\y|^2}d\y,\\
  \tilde\Omega(f,g)=&\F^{-1}\int \cP_{XL+LX}\frac{\al \hat f(\x-\y)\hat{\bar
      g}(\y)}{|\xi-\eta|^2-|\eta|^2-\alpha|\xi|}d\y.
\end{align*}
The equations after normal form reduction can be written as
\begin{equation} \label{EQNF}
  \begin{split}
    & (i\partial_t+D^2)(u-\Omega(N,u))=(Nu)_{LH+HH+\al L}+\Omega(\al D|u|^2,u)+\Omega(N,Nu), \\
    & (i\partial_t+\al D)(N-D\tilde\Omega(u,u))=\al D|u|^2_{HH+\al L+L\al}+D\tilde
    \Omega(Nu,u)+D\tilde\Omega(u,Nu).
  \end{split}
\end{equation}

\subsection{Function spaces and Strichartz estimates}\label{subsect:fs}
Let $s,l \in \R$ and $1\leq p,q\leq \infty$.  We use $B^s_{p,q},\dot B^s_{p,q}$ to denote the standard Besov space, with norms
\[\norm{f}_{B^s_{p,q}}=\norm{P_{\leq 0}f}_p+\Big(\sum_{k=1}^\infty 2^{ksq}\norm{P_kf}_p^q\Big)^{1/q},\quad \norm{f}_{\dot B^s_{p,q}}=\Big(\sum_{k=-\infty}^\infty 2^{ksq}\norm{P_kf}_p^q\Big)^{1/q},\]
with obvious modifications if $q=\infty$, and we simply write $B^s_p=B^s_{p,2}, \dot B^s_p=\dot B^s_{p,2}$.

For the exponents $s\le l+1$, we use the following resolution spaces
\EQ{ \label{def XY}
 \pt u \in X^s:=C(\R;H^s(\R^4))\cap L^\I(\R;H^s(\R^4))\cap L^2(\R;B^s_4(\R^4)),
 \pr N \in Y^l:=C(\R;H^l(\R^4))\cap L^\I(\R;H^l(\R^4)).}
For any Banach function space $Z$ on $\R^{1+4}$ and any interval $I\subset\R$, the restriction of $Z$ onto $I$ is denoted by $Z(I)$. For example,
\EQ{
 X^s([0,T])=C([0,T];H^s(\R^4))\cap L^2((0,T);B^s_4(\R^4)).}

We will use the following well-known Strichartz estimates for the wave
and the Schr\"odinger equation in dimension $d=4$.
\begin{lem}[Strichartz estimates, see \cite{KT}]
For any $s\in\R$ and any functions $\phi(x),f(t,x)$, we have
\begin{align*}
      \norm{S(t)\phi}_{{L_t^\infty H_x^s}\cap {L_t^2B^s_4}}\les &
      \norm{\phi}_{H^s}\\
      \normb{\int_0^tS(t-s)f(s)ds}_{{L_t^\infty L_x^2}\cap
        {L_t^2B^0_4}}\les & \norm{f}_{{L_t^1L_x^2}+{L_t^2B^0_{4/3}}}.
\end{align*}
\begin{align*}
      \norm{W_\al(t)\phi}_{{L_t^\infty L_x^2}\cap {L_t^2\dot
          B_6^{-5/6}}}\les\ &
      \norm{\phi}_{L^2}\\
      \normb{\int_0^tW_\al(t-s)f(s)ds}_{{L_t^\infty L_x^2}\cap
        {L_t^2\dot B_6^{-5/6}}}\les\ & \norm{f}_{{L_t^1L_x^2} }.
\end{align*}
\end{lem}

\subsection{Multi-linear estimates for quadratic and cubic terms}\label{subsect:qc-est}
Next, we prove multi-linear estimates for the nonlinear terms in
\eqref{EQNF} in the Besov spaces of $x\in\R^4$. For $t$, only
H\"older's inequalities will be used, which needs no explanation. In
the following, we ignore the dependence of constants on $(s,l)$, but
distinguish by $C(K)$ when it is not uniform for $K$. The main tools
are Littlewood-Paley theory and certain Coifman-Meyer type bilinear
Fourier multiplier estimates. Roughly speaking, the multipliers
$\Om$ and $\ti\Om$ act like \EQ{ \label{symbol est}
 \Om(f,g) \sim D^{-1}\LR{D}^{-1}(fg)_{XL},
 \pq \ti\Om(f,g) \sim D^{-1}\LR{D}^{-1}(f\bar{g})_{XL+LX},}
in product estimates in the Besov spaces.
Hence the proof is reduced to usual computation of exponents as in the paraproduct. We only sketch the proof.
\begin{lem}[Quadratic terms] \label{BLT} Let  $K \geq 5$. \\
{\rm(1)} Assume that $s,l\ge 0$. Then for any $N(x)$ and $u(x)$,
\EQ{ \label{quad1}
  \pt \|(Nu)_{LH+\al L}\|_{B^s_{4/3}} \les \|N\|_{H^l}\|u\|_{B^s_4},
  \pr \|(Nu)_{HH}\|_{B^s_{4/3}} \les C(K)\|N\|_{H^l}\|u\|_{B^s_4}.}
{\rm(2)} Assume $0\le l+1\leq 2s$. Then for any $u(x)$ and $v(x)$,
\EQ{ \label{quad2}
   \pt \|D(uv)_{HH}\|_{H^l} \les C(K)\|u\|_{B^s_4}\|v\|_{B^s_4},
   \pr \|D(uv)_{\al L + L\al}\|_{H^l} \les \|u\|_{B^s_4}\|v\|_{B^s_4}.}
\end{lem}

\begin{proof}
The estimates above follow directly from Bony's paraproduct and H\"older's inequality.  For example,
\EQ{
    \norm{P_k(Nu)_{LH}}_{L^{4/3}}\les & \sum_{j=k-2}^{k+2}\norm{(P_{\leq j-K}N)(P_ju)}_{L^{4/3}}
    \les  \sum_{j=k-2}^{k+2}\norm{N}_{L^{2}}\norm{P_ju}_{L^4}.}
Then, we sum up the squares with respect to $k$. The other estimates follow in a similar manner.
This argument loses the summability for $HH$ at the $0$ regularity ($s=l=0$ for (1) and $s=l+1=0$ for (2)), but then we can simply use H\"older in $x$ together with the embedding $B^0_p\subset L^p$ and $L^{p'}\subset B^0_{p'}$ for $2\le p\le\I$.
\end{proof}

Similarly to \cite[Lemma 4.4]{GNW} and \cite[Lemma 4.4]{gnw-kg2}, we will exploit in the proof of local well-posedness and persistence of regularity that the boundary contributions, as well as cubic terms, can be made small by choosing $K\geq 5$ large.
\begin{lem}[Boundary terms]\label{lem:boundary}
There exist $\te_j(s,l)\ge 0$ such that for all $K\ge 5$, and for any $N(x),u(x),v(x)$, we have the following:\\
{\rm(1)} If $l\ge \max(0,s-2)$ and $(s,l)\not=(2,0)$,
\EQ{ \label{it-b0}
      \|\Om(N,u)\|_{H^{s}} \les\ & 2^{-\te_1 K} \|N\|_{H^l}\|u\|_{H^s},
 \pq \te_1>0\text{ for $s<l+2$}.}
{\rm(2)} If $l\le \min(2s-1,s+1)$ and $(s,l)\not=(2,3)$,
\EQ{ \label{it-b1}
      \|D\tilde\Om(u,v)\|_{H^l}\les\ & 2^{-\te_2 K} \|u\|_{H^s}\|v\|_{H^s}, \pq \te_2>0\text{ for $l<s+1$.}}
{\rm(3)} If $l\ge\min(0,s-1)$ and $(s,l)\not=(1,0)$,
\EQ{ \label{it-b2}
    \|\Om(N,u)\|_{B^s_4}\les\  2^{-\te_3 K}\|N\|_{H^l}\|u\|_{B^s_4}, \pq \te_3>0\text{ for $s<l+1$.}}
{\rm(4)} If $l\le \min(2s-1/2,s+3/2)$ and $(s,l)\not=(2,7/2)$,
\EQ{ \label{it-b3}
  \|\LR{D}^l\tilde\Om(u,v)\|_{\dot B^{1/6}_{6}}\les\  2^{-\te_4 K}
    [\|u\|_{B^s_4}\|v\|_{H^s}+\|v\|_{B^s_4}\|u\|_{H^s}],}
where $\te_4>0$ for $l<s+3/2$.
\end{lem}
\begin{proof}
Since they are all straightforward, we prove only \eqref{it-b2}-\eqref{it-b3}, leaving \eqref{it-b0}-\eqref{it-b1} to the reader.
By \cite[Lemma 3.5]{GN} and using \eqref{symbol est} with Bernstein, we have
\EQ{ \label{est Om}
 \pt \|P_k\LR{D}D\Om(P_{k_0}N,P_{k_1}u)\|_{L^p_x} \lec \|P_{k_0}N\|_{L^{p_0}_x}\|P_{k_1}u\|_{L^{p_1}_x}
 \pr \lec 2^{4k_0(1/p_0-1/q_0)+4k_1(1/p_1-1/q_1)}\|P_{k_0}N\|_{L^{q_0}_x}\|P_{k_1}u\|_{L^{q_1}_x},}
for any $k,k_0,k_1\in\Z$ and any $p,p_0,p_1,q_0,q_1\in[1,\I]$ satisfying $1/p=1/p_0+1/p_1$ and $q_j\le p_j$. The same estimate holds for the bilinear operator $\ti\Om$. For the low frequency part, say if $k_1\le k_0-K$, we can replace $P_{k_1}$ with $P_{\le k_1}$.
The above with $(p,p_0,p_1,q_0,q_1)=(4,4,\I,2,4)$ and the $HL$ restriction $|k-k_0|\le 1$ in $\Om$ yields
\EQ{
  \|\Om(N,u)\|_{B^s_4}
 \pt\lec \normb{2^{k^+(s-1-l)}\norm{P_kN}_{H^l}\sum_{k_1\leq k-K} 2^{k_1-k_1^+s}\norm{P_{k_1}u}_{B^s_4}}_{l_k^2},}
where $k^+:=\max(k,0)$, using $P_{\le 0}B^s_p\subset \dot B^0_{p,\I}$ for the lower frequency component.
The summation over $k_1\le k-K$ is bounded by
\EQ{
 k\le K \implies 2^{k-K}, \pq k>K \implies \CAS{2^{(1-s)^+(k-K)} &(s\not=1)\\ k-K &(s=1).}}
This and $\|P_kN\|_{H^l}\in\ell^2_k$ lead to \eqref{it-b2}, with the small factor $2^{-\te_3 K}$ for $s<1$ and for $1\le s<l+1$.
The conditions $l\ge 0$ and $l\ge s-1$ ensure uniform boundedness of the coefficient after the summation, respectively for $s<1$ and for $s>1$, while the endpoint $(s,l)=(1,0)$ is excluded due to the logarithmic growth at $s=1$.
Similarly with $(p,p_0,p_1,q_0,q_1)=(6,6,\I,4,2)$, we have
\EQ{
 \pt\|P_k\jb{D}^l \tilde\Om(u,v)_{HL}\|_{\dot B^{1/6}_{6}}
  \pr\lec
     2^{k^+(l-1-s)-k/2}\sum_{k_1\leq k-K}2^{2k_1-k_1^+s}\|P_k u \|_{B^s_4}\|P_{k_1}v\|_{H^s}.}
Using this and $\|P_ku\|_{B^s_4}\in\ell^2_{k\ge 0}$ lead to \eqref{it-b3}, with the small factor for $s<2$ and for $2\le s<l-3/2$.
\end{proof}

\begin{lem}[Cubic terms] \label{lem:cubic}
There exist $\te_j(s,l)\ge 0$ such that for all $K\ge 5$, and for any $M(x),N(x),u(x),v(x),w(x)$, we have the following:\\
{\rm(1)} If $s\geq 1/2$, then $\te_1>0$ and
\EQ{ \label{cub1}
 \|\Om(D(uv),w)\|_{H^s}\les
     2^{-\te_1 K} [\|u\|_{H^s}\|v\|_{B^{1/2}_4}+\|v\|_{H^s}\|u\|_{B^{1/2}_4}]\|w\|_{B^{1/2}_4}.}
{\rm(2)} If $l\geq 0$, $-l<s\le l+2$, $s\le 2l+1$ and $(s,l)\ne (1,0)$,
\EQ{ \label{cub2}
     \|\Om(M,Nu)\|_{B^s_{4/3}} \les 2^{-\te_2 K}\|M\|_{H^l}\|N\|_{H^l}\|u\|_{B^s_4}, \pq \te_2>0 \text{ for $s<l+2$.}}
{\rm(3)} If $s\ge 1/2$, $-s<l\le s+1$, $l\le 2s$, and $(s,l)\ne (1,2)$,
\EQ{ \label{cub3}
     \|D\tilde\Om(Nu,v)\|_{H^l} + \|D\tilde\Om(v,Nu)\|_{H^l} \les 2^{-\te_3 K} \|N\|_{H^l}\|u\|_{B^s_4}\|v\|_{B^s_4},}
where $\te_3>0$ for $l<s+1$.
\end{lem}
\begin{proof}
For \eqref{cub1}, we can use a standard product inequality for $s\ge 1/2$:
\EQ{
 \|uv\|_{B^s_{8/5}} \lec \|u\|_{H^s}\|v\|_{B^{1/2}_4}+\|v\|_{H^s}\|u\|_{B^{1/2}_4},}
which easily follows using $B^{1/2}_4\subset L^8$, e.g.~by the paraproduct calculus. Putting $f:=uv$, we obtain from \eqref{est Om} with $(p,p_0,p_1,q_0,q_1)=(2,2,\I,8/5,4)$
\EQ{
 \|P_k\Om(Df,w)\|_{H^s}
  \lec 2^{k/2-k^+}\sum_{k_1\le k-K}2^{k_1-k_1^+/2}\|f\|_{B^s_{8/5}}\|w\|_{B^{1/2}_4},}
which leads to \eqref{cub1} with a small factor, in the same way as in the previous lemma.

For \eqref{cub2} and \eqref{cub3}, we can use a standard product inequality:
\EQ{
 \s \le \min(s,l,s+l-1) \implies \|Nu\|_{H^\s} \lec \|N\|_{H^l}\|u\|_{B^s_4},}
which holds for $s+l>0$ unless $s=1$ and $\s=l$.
Putting $g:=Nu$, we obtain from \eqref{est Om} with $(p,p_0,p_1,q_0,q_1)=(4/3,2,4,2,2)$
\EQ{ \label{Om Mg}
 \|P_k\Om(M,g)\|_{B^s_{4/3}}
  \lec \sum_{k_1\le k-K}2^{k^+(s-1-l)-k+k_1-k_1^+\s}\|P_kM\|_{H^l}\|P_{k_1}g\|_{H^\s}.}
First, the low frequency part $k\le 0$ is bounded using Young on $\Z$
\EQ{
 \pt\|P_{\le 0}\Om(M,g)\|_{B^s_{4/3}}
 \pn\lec \|P_k\Om(M,g)\|_{\ell^1_{k\le 0}L^{4/3}_x}
 \pr\lec \|P_kM\|_{\ell^2_{k\le 0}H^l_x} \normb{\sum_{k_1\le k-K}2^{-k+k_1}\|P_{k_1}g\|_{H^\s_x}}_{\ell^2_{k\le 0}}
 \lec 2^{-K}\|M\|_{H^l}\|g\|_{H^\s}.}
For $0<k\le K$, the summation over $k_1$ is bounded by $2^{k(s-l-1)-K}\|P_kM\|_{H^l}\in\ell^2_k$ with the small factor for $s<l+2$.
For $K<k$, it is bounded by
\EQ{
 \CAS{\s<1 \implies 2^{k(s-1-l-\s)}2^{-K(1-\s)} \\ \s>1 \implies 2^{k(s-2-l)}.}}
The case $\s<1$ is fine if $\s=l$ by $s\le 2l+1$, if $\s\le s+1$ by $l\ge 0$, and if $\s=s+l-1$ by $l\ge 0$.
In the critical case $s=1$ for the product inequality, we have $s<2l+1$ and $l>0$ by the exclusion $(s,l)\not=(1,0)$, so that we can choose $\s=l-\e$.
The case $\s>1$ is fine by $s\le l+2$.
Then the only remaining case is $(s,l)=(3,1)$, where we are forced to choose $\s=1$ then we should replace \eqref{Om Mg} for $k>K$ with
\EQ{
 \|P_k\Om(M,g)\|_{B^s_{4/3}} \lec 2^{k(s-2-l)}\|P_kM\|_{H^l}\|P_{\le k-K}g\|_{H^1},}
which is bounded using $\|P_kM\|_{H^l}\in\ell^2_k$.
Thus we obtain \eqref{cub2}.

Similarly, from \eqref{est Om} with $(p,p_0,p_1,q_0,q_1)=(2,2,\I,2,4)$, we have
\EQ{
 \pt\|P_kD\ti\Om(g,v)_{HL}\|_{H^l}+\|P_kD\ti\Om(v,g)_{LH}\|_{H^l}
 \pr\lec \sum_{k_1\le k-K}2^{k^+(l-1-\s)+k_1-k_1^+s}\|P_kg\|_{H^\s}\|P_{k_1}v\|_{B^s_4},
}
for which the low frequencies $k\le K$ are easily bounded using the factor $2^{k_1}$, while for $k>K$ the summation is bounded by
\EQ{ \label{est ghigh}
 \CAS{s<1 \implies 2^{k(l-s-\s)}2^{-K(1-s)},\\
 s=1 \implies 2^{k(l-1-\s)}(k-K),\\ s>1 \implies 2^{k(l-1-\s)}.}}
The case $s<1$ is fine if $\s=s$ by $l\le 2s$, and if $\s=s+l-1$ by $s\ge 1/2$.
The case $s>1$ is fine if $\s=s$ by $l\le s+1$, and obviously if $\s=l$.
The critical case $s=1$ is also fine, as none of the conditions is on the boundary thanks to $(s,l)\not=(1,2)$.

For the other $HL$ interaction, choosing $(p,p_0,p_1,q_0,q_1)=(2,4,4,2,4)$ we have
\EQ{
 \pt\|P_kD\ti\Om(g,v)_{LH}\|_{H^l}+\|P_kD\ti\Om(v,g)_{HL}\|_{H^l}
 \pr\lec \sum_{k_1\le k-K}2^{k^+(l-1-s)+k_1-k_1^+\s}\|P_kv\|_{B^s_4}\|P_{k_1}g\|_{H^\s},
}
which is also easy for $k\le K$. For $k>K$, the summation is bounded by
\EQ{
 \CAS{\s<1 \implies 2^{k(l-\s-s)}2^{-K(1-\s)} \\
 \s>1 \implies 2^{k(l-1-s)}.}}
The case $\s<1$ is the same as the case $s<1$ in \eqref{est ghigh}.
The case $\s>1$ is OK by $l\le s+1$.
When $l=s+1\ge 3/2$, we can choose $\s=\min(s,l,s+l-1)=s\not=1$ thanks to $(s,l)\not=(1,2)$. In the critical case $s=1$, we can choose $\s<\min(s,l,s+l-1)\le 1$ such that $l-s-\s<0$, since $l<2s=2$. This concludes the proof of \eqref{cub3}.
\end{proof}

\section{Small data scattering for $s\le l+1$} \label{sect:scat}
Using the multilinear estimates in the previous section, it is now easy to obtain global well-posedness and scattering for small initial data in $H^s\times H^l$ in the range \eqref{sl range} under $s\le l+1$. In Section 5 we will show that we only need smallness in $H^{1/2}\times L^2$ for all regularities by persistence of regularity argument.
Fix $K=5$.
As in \cite[Section~4]{GN}, for fixed initial data $(u_0,N_0)\in H^s\times H^l$, we define a mapping $(u,N)\mapsto (u',N')=\Phi_{u_0,N_0}(u,N)$ by the right-hand sides of the equations \eqref{eq:intu}-\eqref{eq:intN}.
Then for small initial data $(u_0,N_0)$, we see that $\Phi_{u_0,N_0}$ is a contraction in a small ball around $0$ of $X^s\times Y^l$.
Indeed, from the estimates in the previous section, we obtain
\EQ{
 \pt \|u'\|_{X^s} \lec \|u_0\|_{H^s} + \|N\|_{Y^l}\|u\|_{X^s} + \|u\|_{X^s}^3 + \|N\|_{Y^l}^2\|u\|_{X^s},
 \pr \|N'\|_{Y^l} \lec \|N_0\|_{H^l} + \|u\|_{X^s}^2 + \|N\|_{Y^l}\|u\|_{X^s}^2,}
where we need $s\le l+1$ in using \eqref{it-b2} for $\Om(N,u)$.
By the contraction mapping principle, we have a unique solution in a small ball in $X^s\times Y^l$, and the Lipschitz continuity of the solution map $H^s\times H^l\to X^s\times Y^l$ follows from the standard argument.

Now we derive scattering for $(u,N)$ in $H^s\times H^l$, assuming $(s,l)$ satisfying \eqref{sl range}, $(u,N)\in X^{1/2}\times Y^0$ with small norm and 
the scattering of the transformed variables, namely for
\EQ{
 \Psi(u,N):=(u-\Om(N,u),N-D\tilde\Om(u,u))}
there exist $(u_\pm, N_\pm)\in H^s\times H^l$ with small norm in $H^{1/2}\times L^2$ such that
\EQ{\label{eq:scatPsi}
\Psi(u,N)-(S(t)u_\pm, W_\al(t)N_\pm)\to 0 \IN{H^{s}\times H^{l}}\pq(t\to\pm\I).}
In the current case $s\leq l+1$, the latter assumption \eqref{eq:scatPsi} obviously holds in view of that $(u,N)\in X^s\times Y^l$ and the Strichartz estimate with the global bounds on the nonlinear terms.

The bilinear estimate for the normal form in Lemma \ref{lem:boundary} implies that the above transform $\Psi$ is invertible for small data in $H^{1/2}\times L^2$  and bi-Lipschtiz.  
More precisely, for any $(u',N')\in H^{1/2}\times L^2$, the inverse image $\Psi^{-1}(u',N')$ is the fixed points of the map 
\EQ{
 (u,N)\mapsto \Psi_{u',N'}(u,N):=(u'+\Om(N,u),N'+D\ti\Om(u,u)).}
Lemma \ref{lem:boundary} implies that $\Psi_{u',N'}$ is a contraction in a small ball of $H^{1/2}\times L^2$ if $(u',N')$ is small, hence the unique small $(u,N)\in\Psi^{-1}(u',N')$ is given by the iteration 
\EQ{ \label{iteration inverse}
 (u,N)=\lim_{k\to\I}(\Psi_{u',N'})^k(0,0).} 
By \eqref{eq:scatPsi}, we get
\[(u,N)-\Psi^{-1}(S(t)u_\pm, W_\al(t)N_\pm)\to 0 \IN{H^{1/2}\times L^2}\pq(t\to\pm\I).\]
To show the scattering for $(u,N)$, it suffices to show 
\EQ{\label{eq:inverlinear}
\Psi^{-1}(S(t)u_\pm, W_\al(t)N_\pm) \to (S(t)u_\pm, W_\al(t)N_\pm) \IN{H^{1/2}\times L^2}\pq(t\to\pm\I).}
By the construction of inverse, we get 
\EQ{
(u_\pm^n(t),N_\pm^n(t)) \to \Psi^{-1}(S(t)u_\pm, W_\al(t)N_\pm) \IN{L_t^\infty (H^{1/2}\times L^2)}\pq(n\to \I).}
where $(u_\pm^0,N_\pm^0)=(0,0)$, and for $n=1,2,\cdots,$
\begin{align*}
u_\pm^{n+1}=&S(t)u_\pm+\Omega(N_\pm^n,u_\pm^n),\\
N_\pm^{n+1}=&W_\al(t)N_\pm+D\tilde\Omega(u_\pm^n,u_\pm^n).
\end{align*}
Thus, to show \eqref{eq:inverlinear}, it suffices to show for any $n$ 
\EQ{
(u_\pm^n(t),N_\pm^n(t)) \to (S(t)u_\pm, W_\al(t)N_\pm) \IN{H^{1/2}\times L^2}\pq(t\to \pm\I),}
for which by induction on $n$ and bilinear estimates it suffices to show
\EQ{
 (\Om(N_F,u_F),D\tilde\Om(u_F,u_F))\to 0 \IN{H^{s}\times H^{l}}\pq(t\to\pm\I)}
for all free solutions $(u_F,N_F)$ in $H^{s}\times H^{l}$.
The density argument with the bilinear estimate allows us to restrict to the case $u_F(0),N_F(0)\in C_0^\I(\R^4)$, then the above is almost obvious by the dispersive decay of $S(t)$ and $W_\al(t)$ (we omit the details). 

For higher regularity $(s,l)\not=(1/2,0)$, we do not have smallness in $H^s\times H^l$, so we should replace Lemma \ref{lem:boundary} with the following set of estimates 
\EQ{ \label{Omest Besov}
 \pt\|\Om(N,u)\|_{H^s} \lec \|N\|_{H^l}\|u\|_{B_u},
 \pq\|\Om(N,u)\|_{B_u} \lec \|N\|_{B_N}\|u\|_{B_u}, 
 \pr\|D\ti\Om(u,u)_{HL}\|_{H^l} \lec \|u\|_{H^s}\|u\|_{B_u},}
where the Besov spaces $B_u$ and $B_N$ are defined by
\EQ{
 B_u:=B^{s-\e}_p, \pq B_N:=B^{l-\e}_p, \pq 1/p=1/2-\e/4}
for some small $\e>0$ such that $H^s\times H^l\subset B_u\times B_N$ by the sharp Sobolev embedding. 
\eqref{Omest Besov} implies that $\Psi_{u',N'}$ is a contraction with respect to the equivalent norm 
\EQ{
 \|(u,N)\|_Z:=\|u\|_{H^s}+\|N\|_{H^l}+\de^{-2}\|u\|_{B_u}} 
for $0<\de\ll 1$, on the closed set 
\EQ{ 
 F:=\{(u,N)\in H^s\times H^l ; \|(u,N)\|_Z\le 1/\de, \pq \|N\|_{B_N} \le \de, \pq \|u\|_{B_u} \le \de^3\},}
provided that $2(u',N')\in F$. 
Indeed, \eqref{Omest Besov} yields for any $(u,N)\in F$, 
\EQ{ \label{Om small}
 \|(\Om(N,u),D\ti\Om(u,u))\|_{H^s\times H^l} \lec \de^2,
 \pq \|\Om(N,u)\|_{B_u} \lec \de^4,}
hence $\|(\Om(N,u),D\ti\Om(u,u))\|_Z \lec \de^2$ and $\Psi_{u',N'}(u,N)\in F$. For the difference, we have from \eqref{Omest Besov}, for any $(v,M)\in H^s\times H^l$, 
\EQ{
 \pt\|(\Om(N,v),D\ti\Om(u,v)_{HL})\|_{H^s\times H^l}
   \lec \|(u,N)\|_{H^s\times H^l}\|v\|_{B_u} \lec \de\|(v,M)\|_{Z},
 \pr\|(\Om(M,u),D\ti\Om(v,u)_{HL})\|_{H^s\times H^l}
   \lec \|(v,M)\|_{H^s\times H^l}\|u\|_{B_u} \lec \de^3\|(v,M)\|_{Z},
 \pr\|\Om(N,v)+\Om(M,u)\|_{B_u}
   \lec \|N\|_{B_N} \|v\|_{B_u} + \|M\|_{B_N} \|u\|_{B_u} \lec \de^3 \|(v,M)\|_{Z}. }
Since the scattering of $\Psi(u,N)$ implies $\|\Psi(u,N)\|_{B_u\times B_N}\to 0$ as $t\to\I$, choosing $\de>0$ small enough ensures that $2\Psi(u,N)\in F$ for large $t$. 
Then $(u,N)$ given by \eqref{iteration inverse} is the same as the fixed point in $F$. 
Since we can take $\de>0$ arbitrarily small, \eqref{Om small} implies that $\|(u,N)-\Psi(u,N)\|_{H^s\times H^l}\to 0$ as $t\to\I$, hence the scattering of $(u,N)$ in $H^s\times H^l$. 

Since all the estimates are uniform and global in time, the same argument works for the final state problem, namely to find the solution for a prescribed (small) scattering data at $t=\I$.
Thus we obtain
\begin{prop} \label{small data scatter}
Let $(s,l)\in\R^2$ satisfy \eqref{sl range}, $s\le l+1$ and $(s,l)\not=(1,0)$. Then there exists $\e_1=\e_1(s,l)>0$ such that for any $(u_0,N_0)\in H^s(\R^4)\times H^l(\R^4)$ satisfying $\|(u_0,N_0)\|_{H^s\times H^l}\le\e_1$, there exists a unique global solution $(u,N)\in X^s\times Y^l$ of \eqref{eq:Zak1}. Moreover, there exists $(u^+,N^+)\in H^s\times H^l$ such that
\EQ{ \label{scattering}
 \lim_{t\to\I}\|u(t)-S(t)u^+\|_{H^s_x}+\|N(t)-W_\al(t)N^+\|_{H^l_x}=0.}
Conversely, for any $(u^+,N^+)\in H^s\times H^l$ with $\|(u^+,N^+)\|_{H^s\times H^l}\le\e_1$, there exists a unique solution $(u,N)\in X^s\times Y^l$ satisfying \eqref{scattering}. Both the maps $(u_0,N_0) \mapsto (u,N)$ and $(u^+,N^+)\mapsto (u,N)$ are Lipschitz continuous from the $\e_1$-ball into $X^s\times Y^l$.
\end{prop}
The uniqueness without the smallness is proved in the next section.
For the question if $(u,N)$ obtained above really solves the equation \eqref{eq:Zak1} before the normal form, see Remark \ref{original eq}.

\section{Large data local well-posedness for $s<l+1$}\label{sect:lwp-large}
For large data, the proof in the previous section does not immediately work, in particular at the endpoint $(s,l)=(1/2,0)$.
The main difficulty is the lack of flexibility in the choice of the Strichartz norm for the boundary term and the bilinear term $(Nu)_{LH}$.
More precisely, $L^\I_tH^l_x$ is the only choice among the Strichartz norms of $W_\al(t)$ for $N$, to estimate $\Om(N,u)$ in $L^\I_t H^s_x$, and to avoid losing regularity in $(Nu)_{LH}$.
For the former term, we can play with the frequency gap parameter $K$ in the normal form to extract a small factor.
For the latter term, we use the following
\begin{lem} \label{N small decop}
Let $0<T\le \I$ and $N\in C([0,T);L^2(\R^4))$. Suppose that $W_\al(-t)N(t)$ is strongly convergent in $L^2_x$ as $t\to T-0$. Then for any $\e>0$, there exists a finite increasing sequence $0=T_0<T_1<\cdots<T_{n+1}=T$ such that
\EQ{
 \|N\|_{(L^\I_t L^2_x+L^2_t L^4_x)(T_j,T_{j+1})}<\e}
for each $j=0\etc n$.
\end{lem}
Note that the $L^2_tL^4_x$ norm is not controlled by the Strichartz estimate for $W_\al(t)$, but it is bounded for nice initial data.
The case $T=\I$ will be used for large data scattering.
For $T<\I$, the assumption on $N$ is equivalent to $N\in C([0,T];L^2)$.
\begin{proof}
Put $N^+:=\lim_{t\to T-0}W_\al(-t)N(t)\in L^2_x$. By the strong convergence, there exists $T'\in(0,T)$ such that
$\sup_{T'\le t<T}\|N(t)-W_\al(t)N^+\|_{L^2_x}<\e/4$.
Since $C_0^\I\subset L^2_x$ is dense, there exists $N_0\in C_0^\I$ such that $\|N_0-N^+\|_{L^2_x}<\e/4$. The dispersive decay of $W_\al(t)$ implies that $W_\al(t)N_0\in L^2_t L^4_x(\R)$. Define $N'$ by
\EQ{
 N'(t):= P_{\le k}N(t)\ (0\le t\le T'), \pq N'(t):= W_\al(t)N_0\ (T'<t<T).}
By the above choice of $T'$ and $N_0$, we have
$\|N-N'\|_{L^\I((T',T);L^2_x)}< \e/2$.
Since $N\in C([0,T'];L^2_x)$ and $[0,T']$ is compact, $N'(t)\to N(t)$ in $L^2_x$ uniformly on $t\in[0,T']$ as $k\to\I$. Hence for large $k$ we have $\|N-N'\|_{L^\I([0,T'];L^2_x)}<\e/2$.
Hence
\EQ{
 \|N-N'\|_{L^\I([0,T);L^2_x)}<\e/2, \pq N'\in L^2_t([0,T);L^4_x).}
Choosing $T_1<T_2<\cdots <T_n$ appropriately ensures that
$\|N'\|_{L^2_t L^4_x(T_j,T_{j+1})}<\e/2$
for each $j$, then we get the desired estimate.
\end{proof}

Now we are ready to prove the local well-posedness for large data in $H^{1/2}\times L^2$.
For any initial data $(u_0,N_0)\in H^{1/2}\times L^2$, let
\EQ{
 u_F:=S(t)(u_0-\Om(N_0,u_0)), \pq N_F:=W_\al(t)(N_0-D\ti\Om(u_0,u_0)),}
and apply Lemma \ref{N small decop} to $N_F$. Then for any $\e>0$, there exists $T>0$ such that
\EQ{ \label{LWP small}
 \|u_F\|_{L^2_tB^{1/2}_4(0,T)}+\|N_F\|_{L^\I_tL^2_x+L^2_tL^4_x(0,T)}<\e.}
Putting $\cH:=H^{1/2}\times L^2$ and $m:=\|(u_F(0),N_F(0))\|_{\cH}$, we look for a unique local solution on $(0,T)$ as a fixed point of the map $\Phi_{u_0,N_0}$ in the closed set
\EQN{
 K_m^\e:=\Biggl\{(u,N)\in C([0,T];\cH) \Biggm|
 \begin{split}
 \pt \|(u,N)\|_{L^\I_t(0,T;\cH)}\le 2m,
 \pr \|u\|_{L^2_tB^{1/2}_4(0,T)}+\|N\|_{L^\I_tL^2_x+L^2_tL^4_x(0,T)}\le 2\e
 \end{split} \Biggr\}.}
From the multilinear estimates in Section \ref{sect:est}, we have
\EQ{
 \pt\|\Om(N,u)\|_{X^{1/2}} \lec 2^{-\te K}\|N\|_{L^\I_tL^2_x}\|u\|_{X^{1/2}},
 \pr\|\Om(D|u|^2,u)\|_{L^1_tH^{1/2}_x} \lec 2^{-\te K}\|u\|_{L^\I_tH^{1/2}}\|u\|_{L^2_tB^{1/2}_4}^2,
 \pr\|\Om(N,Nu)\|_{L^2_tB^{1/2}_{4/3}} \lec 2^{-\te K}\|N\|_{L^\I_tL^2_x}^2\|u\|_{L^2_tB^{1/2}_4},
 \pr\|D\ti\Om(u,u)\|_{L^\I_tL^2_x}\lec 2^{-\te K}\|u\|_{L^\I_tH^{1/2}_x}^2,
 \pr\|D\ti\Om(Nu,u)\|_{L^1_tL^2_x}\lec 2^{-\te K}\|N\|_{L^\I_tL^2_x}\|u\|_{L^2_tB^{1/2}_4}^2,}
and the same estimate on $D\ti\Om(u,Nu)$, as well as for the difference.
Taking $K$ large makes these estimates contractive.
For the remaining two terms,
\EQ{
 \pt\|D|u|^2_{HH+\al L+L\al}\|_{L^1_tL^2_x}\lec C(K)\|u\|_{L^2_tB^{1/2}_4}^2,
 \pr\|(Nu)_{LH+HH+\al L}\|_{L^2_tB^{1/2}_{4/3}+L^1_tH^{1/2}_x} \lec C(K)\|N\|_{L^\I_tL^2_x+L^2_tL^4_x}\|u\|_{L^2_tB^{1/2}_4},}
which is also made contractive on the interval $[0,T]$ by choosing $\e>0$ small enough such that $C(K)\e\ll 1$ after fixing $K$.
Then $\Phi_{u_0,N_0}$ becomes a contraction on $K_m^\e$.

The uniqueness of solution in the class $X^{1/2}\times Y^0$ is obtained in the same fashion: Let $(u_j,N_j)\in X^{1/2}\times Y^0$ for $j=0,1$ be two solutions. For any $\e>0$, applying Lemma \ref{N small decop}, we can find $T'\in(0,T)$ such that for $j=0,1$
\EQ{
 \|u_j\|_{L^2_tB^{1/2}_4(0,T')}+\|N_j\|_{L^\I_tL^2_x+L^2_tL^4_x(0,T')}<\e,}
so that both the solutions belong to $K_m^\e$ on $[0,T']$, hence $(u_0,N_0)=(u_1,N_1)$ as long as they are solutions in the above class.

The continuous dependence is also obtained in the same way, because
\EQ{
 H^{1/2}\times L^2 \ni (u_0,N_0) \mapsto (u_F,N_F)\in X^{1/2}\times Y^0}
is continuous. Take a strongly convergent sequence of initial data. If the smallness condition \eqref{LWP small} is satisfied by the limit, then so is it by those sufficiently close to the limit. Then we can estimate the difference from the limit in the same way as above, leading to the strong continuity.

We have worked out at the lowest regularity $(s,l)=(1/2,0)$, but the same argument works as long as we have the small factor $2^{-\te K}$, namely for $|s-l|<1$. Thus we obtain
\begin{prop} \label{LWP large}
Let $(s,l)\in\R^2$ satisfy \eqref{sl range} and $|s-l|<1$.
For any $(u_0,N_0)\in H^s(\R^4)\times H^l(\R^4)$, there exists a unique local solution $(u,N)\in (X^s\times Y^l)([0,T])$ of \eqref{eq:Zak1} for some $T>0$, where both $T$ and $(u,N)$ depend continuously on $(u_0,N_0)$. More precisely, if $(u_{0,n},N_{0,n})\to (u_0,N_0)$ in $H^s\times H^l$, then $T_n\to T$ and for any $0<T'<T$, we have
$\|u_n-u\|_{X^s([0,T'])}+\|N_n-N\|_{Y^l([0,T'])}\to 0$.
\end{prop}

\section{Persistence of regularity except for $(s,l)=(1,0)$}
Once we have the unique solution at the lowest regularity $(s,l)=(1/2,0)$, it gains as much regularity as the initial data.
To prove this, we will focus on the derivation of a priori estimates, assuming that all relevant norms are finite, which is justified by the local well-posedness in higher regularity by Proposition \ref{LWP large}.

For solutions $(u,N)\in(X^{1/2}\times Y^0)([0,T))$ with $(u(0),N(0))\in H^s\times H^l$ and $0<T\le\I$, we will improve the regularity up to $H^s\times H^l$ by the following steps.
\begin{enumerate}
\item Improve $u$ to $s<l+1$.
\item Improve $N$ to $l\le 2s-1$, $l\le s+1$, and $(s,l)\not=(2,3)$, for $s<l+1$.
\item Improve $u$ to $1<s<4l+1$, $s\le 2l+11/8$ and $s\le l+2$.
\end{enumerate}
The persistence of regularity is a general phenomenon in nonlinear wave equations, but we encounter some difficulties.
One is the same as in the previous section, which is solved by Lemma \ref{N small decop}.
Another difficulty for $s\ge l+1$ is that the normal form can not keep the full Strichartz norm of $u$, which is why we separate (3).

\subsection{Regularity upgrade for $u$ in $s<l+1$} \label{ss:reg1}
Let $(s,l)\in\R^2$ satisfy \eqref{sl range} and $s<l+1$.
Let $(u_0,N_0)\in H^s\times H^l$ and let $(u,N)\in(X^{1/2}\times Y^l)([0,T))$ be a solution for some $0<T\le 0$.
If $T=\I$, we also assume that $N$ scatters in $H^l_x$.
From the estimates in Section \ref{sect:est}, we have for $s<l+1$,
\EQ{
 \pt \|(Nu)_{LH+HH+\al L}\|_{L^2_tB^s_{4/3}+L^1_tH^s_x} \le C_1(K)\|N\|_{L^\I_t L^2_x+L^2_tL^4_x}\|u\|_{L^2_tB^s_4},
 \pr \|\Om(N,u)\|_{L^\I_tH^s_x} \le C_02^{-\te K}\|N\|_{L^\I_tL^2_x}\|u\|_{L^\I_t H^s_x},
 \pr \|\Om(N,u)\|_{L^2_t B^s_4} \le C_02^{-\te K}\|N\|_{L^\I_tH^l_x}\|u\|_{L^2_t B^s_4},
 \pr \|\Om(D|u|^2,u)\|_{L^1_tH^s_x} \le C_02^{-\te K}\|u\|_{L^\I_t H^s_x}\|u\|_{L^2_t B^{1/2}_4}^2,
 \pr \|\Om(N,Nu)\|_{L^2_tB^s_{4/3}} \le C_02^{-\te K}\|N\|_{L^\I_tH^l_x}^2\|u\|_{L^2_t B^s_4},}
for some constants $\te(s,l)>0$, $C_0(s,l)>0$ and $C_1(K,s,l)>0$. Note that $C_0\to\I$ as $(s,l)\to(1,0)$ in the third and the last estimates, and the small factor $2^{-\te K}$ is lost for $s=l+1$ in the third estimate.
Anyway, taking $K=K(s,l)$ large ensures smallness of the right side in the latter 4 estimates:
\EQ{ \label{small for regup}
 C_02^{-\te K}\BR{\|N\|_{L^\I_tH^l_x}+\|u\|_{L^2_t B^{1/2}_4}^2+\|N\|_{L^\I_tH^l_x}^2} \ll 1.}
After fixing such $K$, choose $\e>0$ such that $C_1(K)\e \ll 1$,
and apply Lemma \ref{N small decop} to $N$, which yields a finite sequence $0=T_0<T_1<\cdots<T_{n+1}=T$ such that
\EQ{ \label{small interval 1}
 \|N\|_{(L^\I_tL^2_x+L^2_tL^4_x)(T_j,T_{j+1})}<\e.}
Then on each subinterval we obtain from the above estimates
\EQ{
 \|u\|_{X^s(T_j,T_{j+1})} \le C_2\|u(T_j)\|_{H^s} + \|u\|_{X^s(T_j,T_{j+1})}/2,}
for some constant $C_2(s)>0$. Hence if $u(0)\in H^s$, then by induction on $j$, we deduce that $u\in X^s([0,T))$. If $T=\I$, this implies the scattering of $u$ in $H^s$, via the argument in Section \ref{sect:scat}.

For continuous dependence on the initial data, consider a sequence of solutions $(u_n,N_n)$ such that $(u_n(0),N_n(0))\to(u(0),N(0))$ in $H^s\times H^l$, $u_n\to u$ in $X^{1/2}(I)$ and $N_n\to N$ in $Y^l(I)$ for some interval $I\subset[0,T)$.
For large $n$, $(u_n,N_n)$ satisfies similar bounds to \eqref{small for regup} and \eqref{small interval 1} within $I$, with slightly bigger bounds.
Then the same estimates as above for $(u_n-u,N_n-N)$ yield the convergence in $(X^s\times Y^l)(I)$.

\subsection{Regularity upgrade for $N$ in $s<l+1$} \label{ss:reg2}
Let $(s,l)\in\R^2$ satisfy \eqref{sl range} and $s<l+1$.
Let $(u_0,N_0)\in H^s\times H^l$ and let $(u,N)\in (X^s\times Y^{l'})([0,T))$ for some $0<T\le \infty $ and some $l'\in(s-1,l)$.
From the estimates in Section \ref{sect:est}, we have
\EQ{
 \pt \|D(|u|^2)_{HH+L\al+\al L}\|_{L^1_tH^l_x} \le C_1(K)\|u\|_{L^2_t B^s_4}^2,
 \pr \|D\ti\Om(u,u)\|_{L^\I_tH^l_x} \le C_0\|u\|_{L^\I_tH^s_x}^2,
 \pr \|D\ti\Om(Nu,u)\|_{L^1_tH^l_x}+\|D\ti\Om(u,Nu)\|_{L^1_tH^l_x} \le C_0\|N\|_{L^\I_tH^l_x}\|u\|_{L^2_tB^s_4}^2,}
for some constants $C_0(s,l)>0$ and $C_1(K,s,l)>0$, and the same for $D\ti\Om(u,Nu)$. Choose $\e>0$ so small that $C_0 \e^2 \ll 1$.
Since $u\in L^2_tB^s_4(0,T)$, there exists a finite sequence $0=T_0<T_1<\cdots<T_{n+1}=T$ such that
\EQ{ \label{small interval 2}
 \|u\|_{L^2_tB^s_4(T_j,T_{j+1})}<\e}
for each $j$. Then on each subinterval we have from the above estimates
\EQ{
 \|N\|_{L^\I_tH^l_x(T_j,T_{j+1})} \pt\le C_2\|N(T_j)\|_{H^l}+\|N\|_{L^\I_tH^l_x(T_j,T_{j+1})}/2
 \prQ+C_1(K)\e^2+C_0\|u\|_{L^\I_tH^s_x(T_j,T_{j+1})}^2,}
for some constant $C_2(l)>0$. Hence if $N(0)\in H^l$, then by induction on $j$, we deduce that $N\in L^\I_tH^l(0,T)$.
If $T=\I$, then we have the scattering of $N$ from the argument in Section \ref{sect:scat}. We also obtain the Strichartz norm of $N$ using \eqref{it-b3} for the normal form.
We can also upgrade continuous dependence, using the same estimates for the difference from the limit, see the previous subsection for more detail.
Combining the results in this and the previous subsections yields
\begin{prop} \label{regup Xs}
Let $(s,l)\in\R^2$ satisfy \eqref{sl range} and $s<l+1$. Let $(u,N)\in(X^{1/2}\times Y^0)(I)$ be a solution of \eqref{eq:Zak1} on an interval $I\subset\R$, and suppose that $(u(t_0),N(t_0))\in H^s(\R^4)\times H^l(\R^4)$ at some $t_0\in I$. Then $(u,N)\in (X^s\times Y^l)(I)$ and moreover,
\EQ{ \label{Stz on N}
  N \in L^2_t(I;\dot B^{l-5/6}_6\cap \dot B^{-5/6}_6).}
If $I\supset(t_0,\I)$, then $(u,N)$ scatters in $H^s\times H^l$ as $t\to\I$.
If $(u_n(t_0),N_n(t_0))\to (u(t_0),N(t_0))$ in $H^s\times H^l$ and the corresponding sequence of solutions $(u_n,N_n)\to(u,N)$ in $(X^{1/2}\times Y^0)(J)$ on some interval $t_0\in J\subset I$, then the convergence holds in $(X^s\times Y^l)(J)$. The same convergence result holds for the scattering data, if $I\cap J\supset(t_1,\I)$ for some $t_1<\I$.
\end{prop}

\subsection{Regularity upgrade for $u$ in $s\ge l+1$} \label{ss:reg3}
Let $(s,l)\in\R^2$ satisfy $s\ge l+1$. Then $l>0$ and $s>1$.
Let $(u_0,N_0)\in H^s\times H^l$ and let $(u,N)\in (X^{s'}\times Y^l)([0,T))$ for some $0<T\le \infty$ and some $s'\in(1,s)$.
In this case, the normal form estimate is not good enough to keep the full Strichartz bound of $u$. Hence we decompose
\EQ{ \label{eq uNu'}
 \CAS{u = u' + \Om(N,u),\\
  (i\p_t-\De)u' = (Nu)_{LH\cdots} + \Om(\al D|u|^2,u) + \Om(N,Nu),} }
where $LH\cdots:=LH+HH+\al L$ for brevity, and look for closed estimates in
\EQ{ \label{part Str}
 \pt u' \in X^s, \pq u \in X':=L^\I_t H^s_x \cap L^{2/(1-\ga)}_t L^\I_x,
 \pr N \in L^\I_t H^l_x \cap L^{2/\ga}_t\B, \pq \B:=\dot B^{l-5\ga/6}_{q_1}\cap \dot B^{-5\ga/6}_{q_1},}
where $1/q_1:=1/2-\ga/3$, for some $\ga\in[0,3/4]$ satisfying
\EQ{ \label{al range}
 \ga+1<s, \pq 2l + \ga/2 + 1 \ge s.}
Such $\ga$ exists if and only if $1<s<4l+1$ and $s\le 2l+11/8$.
Also note that
\EQ{
 X^s \subset L^{2/(1-\ga)}_t B^s_{4/(1+\ga)} \subset L^{2/(1-\ga)}_t L^\I_x}
since $\ga+1<s$.
Similarly, $L^{2/\ga}_t\B$ is a wave-Strichartz norm in $H^l$, cf. \eqref{Stz on N}.

We write $(Nu)_{LH\cdots}= (Nu')_{LH\cdots}+ (N\Om(N,u))_{LH\cdots}$. From the estimates in Section \ref{sect:est}, we have for $k>2$,
\EQ{
 \pt\|(Nu')_{LH\cdots}\|_{L^2_tB^s_{4/3}+L^1_tH^s_x} \le C_1(K)\|N\|_{L^\I_t L^2_x+L^2_tL^4_x}\|u'\|_{L^2_tB^s_4},
 \pr\|\Om(D|u|^2,u)\|_{L^1_tH^s_x} \le C_02^{-\te K}\|u\|_{L^2_tB^{1/2}_4}^2\|u\|_{L^\I_tH^s_x},
 \pr\|P_{>k}\Om(N,u)\|_{L^\I_tH^s_x} \le C_0\|N_{>k-1}\|_{L^\I_tH^l_x}\|u\|_{L^\I_t H^s_x},
}
for some constants $C_0(s,l)>0$, $\te(s)>0$ and $C_1(K,s,l)>0$.
We need some more estimates. Since $H^{l+2}\subset L^\I_x$, we have for $k>2$,
\EQ{ \label{part Str Om}
 \pt\|P_{>k}\Om(N,u)\|_{L^{2/(1-\ga)}_t L^\I_x} \le C_0\|N_{>k-1}\|_{L^\I_tH^l_x}\|u\|_{L^{2/(1-\ga)}_t L^\I_x}.}
It remains to estimate $\Om(N,Nu),(N\Om(N,u))_{LH\cdots}$.
If $\B\subset L^4_x$, then for $k\gg\LR{\log\al}$,
\EQ{
 \pt\|P_{>k}\Om(N,Nu)\|_{B^{l+2}_{4/3}} \lec \|N_{>k-1}\|_{H^l_x}\|Nu\|_{L^4_x}
  \lec \|N_{>k-1}\|_{H^l_x}\|N\|_{\B}\|u\|_{L^\I_x},
 \pr\|P_{>k}(N\Om(N,u))_{LH\cdots}\|_{B^{l+2}_{4/3}} \lec \|N\|_{\B}\|N_{>k-K-3}\|_{H^l_x}\|u\|_{L^\I_x}.}
If $\B\not\subset L^4_x$ but $l\ge 5\ga/6$, then putting
\EQ{
 1/q_2:=1/q_1-(l-5\ga/6)/4=1/2-\ga/8-l/4,
 \pq 1/q_3:=1/2+1/q_2,}
we have $\B\subset L^{q_2}$, $B^s_{4/3}\supset B^{2l+\ga/2+1}_{4/3}\supset B^{l+2}_{q_3}$, and so
\EQ{
 \pt \|P_{>k}\Om(N,Nu)\|_{B^s_{4/3}} \lec \|N_{>k-1}\|_{H^l_x}\|Nu\|_{L^{q_2}}
 \lec \|N_{>k-1}\|_{H^l_x}\|N\|_{\B}\|u\|_{L^\I_x},
 \pr \|P_{>k}(N\Om(N,u))_{LH\cdots}\|_{B^s_{4/3}} \lec \|N\|_{\B}\|N_{>k-K-3}\|_{H^l}\|u\|_{L^\I_x}. }
If $l<5\ga/6$, then using
\EQ{
 \pt B^s_{4/3}\supset B^{2l+\ga/2+1}_{4/3}\supset B^{2l-5\ga/6+2}_{q_4},
  \pq 1/q_4:=1/2+1/q_1=1-\ga/3,}
we have
\EQ{
 \|P_{>k}(N\Om(N,u))_{LH\cdots}\|_{B^s_{4/3}} \pt\lec \|N\|_{\B}\|\Om(N,u)_{>k-K-2}\|_{H^{l+2}}
 \pr\lec \|N\|_{\B}\|N_{>k-K-3}\|_{H^l}\|u\|_{L^\I_x}.}
For the other term, putting $\s:=5\ga/6-l>0$ and $\be:=l/(l+\s)\in(0,1)$,
we have the complex interpolation
\EQ{
 \pt [H^l,B^{-\s}_{q_1}]_\be = B^0_{q_5}\subset L^{q_5},
   \pq  [H^l,B^{-\s}_{q_1}]_{1-\be} = B^{l-\s}_{q_6},}
where $1/q_5:=(1-\be)/2+\be/q_1$ and $1/q_6:=1/q_4-1/q_5$, whereas
\EQ{
 \pt\|P_{>k}\Om(N,Nu)\|_{B^s_{4/3}} \pn\lec \|P_{>k}\Om(N,Nu)\|_{B^{l-\s+2}_{q_4}}
  \pr\lec \|N_{>k-1}\|_{B^{l-\s}_{q_6}}\|Nu\|_{L^{q_5}}
  \lec \|N_{>k-1}\|_{B^{l-\s}_{q_6}}\|N\|_{B^0_{q_5}}\|u\|_{L^\I_x}.}
Hence by the interpolation inequality,
\EQ{
 \|\Om(N,Nu)\|_{B^s_{4/3}} \lec \|N_{>k-1}\|_{H^l}^{\be}\|N\|_{H^l}^{1-\be}\|N\|_{\B}\|u\|_{L^\I_x}.}
Therefore, in any case we have some $\be(l,\ga)\in(0,1]$ such that
\EQ{
 \pt\|P_{>k}\Om(N,Nu)\|_{L^2_tB^s_{4/3}}
  \le C_2\|N_{>k-1}\|_{L^\I_tH^l_x}^{\be}\|N\|_{L^\I_tH^l_x}^{1-\be}\|N\|_{L^{2/\ga}_t\B}\|u\|_{X'},
 \pr\|P_{>k}(N\Om(N,u))_{LH\cdots}\|_{L^2_tB^s_{4/3}}
  \le C_2\|N_{>k-K-3}\|_{L^\I_tH^l_x}\|N\|_{L^{2/\ga}_t\B}\|u\|_{X'},}
for some constant $C_2(s,l)>0$.
Choose $K\gg 1$ so large that $C_02^{-\te K}\|u\|_{L^2_tB^{1/2}_4}^2 \ll 1$, and then choose $\e>0$ so small and $k\gg K$ so large that
\EQ{ \label{small interval 3}
 \pt C_1(K)\e + C_0\|N_{>k-1}\|_{L^\I_tH^l_x} \ll 1,
 \pr C_2\|N_{>k-K-3}\|_{L^\I_tH^l_x}^{\be}\|N\|_{L^\I_tH^l_x}^{1-\be}\|N\|_{L^{2/\ga}_t\B} \ll 1.}
Applying Lemma \ref{N small decop}, we obtain a finite sequence $0=T_0<T_1<\cdots <T_{n+1}=T$ such that \eqref{small interval 1} holds.
Then from the above estimates on each subinterval,
\EQ{
 \pt\|u'_{>k}\|_{X^s(T_j,T_{j+1})} \le C_3\|u'(T_j)\|_{H^s} + \de\|u'\|_{X^s(T_j,T_{j+1})} + \de\|u\|_{X'(T_j,T_{j+1})},
 \pr\|\Om(N,u)_{>k}\|_{X'(T_j,T_{j+1})}  \le \de\|u\|_{X'(T_j,T_{j+1})}, }
for some small constant $\de>0$, while the frequencies $\le k$ are bounded by $X^{1/2}$. Using $u=u'+\Om(N,u)$ and $X^s\subset X'$, and adding the low frequencies, we obtain
\EQ{
 \pt\|u'\|_{X^s(T_j,T_{j+1})}+\|\Om(N,u)\|_{X'(T_j,T_{j+1})}
  \prQ\le 2C_3\|u'(T_j)\|_{H^s} + 2^{k(s-1/2)}\|u\|_{X^{1/2}(T_j,T_{j+1})}.}
By induction on $j$ starting from $\|u'(0)\|_{H^s}<\I$, we thus obtain
\EQ{
 \|u\|_{X'(0,T)} \lec \|u'\|_{X^s(0,T)} + \|\Om(N,u)\|_{X'(0,T)} <\I.}
If $T=\I$, then $u$ is scattering by the argument in Section \ref{sect:scat}.
Thus we have obtained
\begin{prop} \label{regup higher s}
Let $(s,l)\in\R^2$ satisfy \eqref{sl range} and $s\ge l+1$. Let $(u,N)\in(X^{1/2}\times Y^0)(I)$ be a solution of \eqref{eq:Zak1} on an interval $I\subset\R$, and suppose that $(u(t_0),N(t_0))\in H^s(\R^4)\times H^l(\R^4)$ at some $t_0\in I$. Then we have $u-\Om(N,u) \in X^s(I)$, as well as \eqref{Stz on N}, and for all $\ga\in[0,3/4]$ satisfying \eqref{al range},
\EQ{ \label{high s space of u}
 \pt u \in C(I;H^s_x)\cap L^\I_t(I;H^s_x)\cap L^{2/(1-\ga)}(I;L^\I_x).}
We also have scattering and continuous dependence similar to Proposition \ref{regup Xs}, but in the above space \eqref{high s space of u}.
\end{prop}
It is easy to replace $L^\I_x$ with $B^s_{4/(1+\ga)}+H^{l+2}$ using $u'\in X^s$ and \eqref{part Str Om}.

\subsection{Lipschitz continuity of the solution map}
Here we consider local Lipschitz continuity of the flow map.
In the above arguments, the Lipschitz dependence is lost only when we seek time intervals with smallness, typically by Lemma \ref{N small decop}.
If $(u_0,N_0)\in H^s\times H^l$ with $s>1/2$ and $l>0$, however, it is easy to see that \eqref{LWP small} holds locally uniformly with respect to the initial data, because we can first dispose of the high frequencies using the higher regularity, and then the remaining low frequencies by Sobolev in $x$ and H\"older in $t$.

Similarly, the regularity upgrading argument in Section \ref{ss:reg1} works uniformly if $l>0$ and $T<\I$, because of \eqref{small interval 1}, so does the argument in Section \ref{ss:reg2} for $s>1/2$, $l<\min(2s-1,s+1)$, and $T<\I$, because of \eqref{small interval 2}, as well as that in Section \ref{ss:reg3} for $l>0$, $s<\min(2l+11/8,l+2)$, and $T<\I$, because of \eqref{small interval 3} and \eqref{small interval 1}.

Thus we obtain Lipschitz continuity of the flow map, locally both in time and in the initial data, for all the exponents $(s,l)$ in the range and off the boundary.
Since we need to decrease $l$ for the uniform control in \eqref{small interval 3}, $\ga$ in \eqref{part Str} can not be on the boundary, namely $2l+\ga/2+1>s$, for the local Lipschitz estimate.

The Lipschitz continuity global in time and for the scattering is more tricky, because the $L^2_tL^4_x$ norm in Lemma \ref{N small decop} is not bounded by the Strichartz estimate for the wave equation.
For small data, we can obtain Lipschitz estimates directly from the contraction mapping argument, but then the smallness on $H^{1/2}\times L^2$ depends on $(s,l)$, which tends to $0$ as $(s,l)$ approaches $s=4l+1$, $(2,3)$ or $(\I,\I)$.
The regularity upgrading for $N$ in Section \ref{ss:reg2} works well for $T=\I$, because in \eqref{small interval 2} the number of subintervals can be uniformly bounded for each $\e>0$, provided that $\|u\|_{L^2_tB^s_4}$ is uniformly bounded. This yields a smallness condition in the form
\EQ{
 \|(u_0,N_0)\|_{H^{1/2}\times L^2} \le \e_2(s,l),}
where $\e_2(s,l)>0$ is non-decreasing in $l$, for global Lipschitz continuity in $H^s\times H^l$.

\begin{rem} \label{original eq}
Strictly speaking, we need to prove that the solution to \eqref{EQNF} obtained above is also a solution of the equation \eqref{eq:Zak1} before the normal form. The easiest way is to use \cite{GTV} for existence of solutions for smooth approximating initial data, taking the limit by the continuous dependence proved above. To be self-contained, however, we can directly show that smooth solutions of \eqref{EQNF} solve \eqref{eq:Zak1}. In fact, if $(u,N)\in(X^s\times Y^s)(I)$ with $s\gg 1$ is a solution of \eqref{EQNF} on some interval $I$, then by definition of $\Om$ and $\ti\Om$, \eqref{EQNF} reads
\EQN{
 \pt eq_u:=(i\p_t+D^2)u-Nu =  -\Om(eq_N,u)-\Om(N,eq_u),
 \pr eq_N:=(i\p_t+\al D)N-\al D|u|^2 =  - D\ti\Om(eq_u,u) - D\ti\Om(u,eq_u).}
Since $eq_u,eq_N\in C(I;H^{s-2})$ and $\Om,D\ti\Om:(H^{s-2})^2\to H^{s-2}$ has a small factor due to $K$, we deduce that $eq_u=0=eq_N$ on $I$ if $K$ is large enough.
\end{rem}

\section{Small data scattering in the energy space}\label{sect:energ}
For $(s,l)=(1,0)$, the failure of Strichartz bound on the normal form $\Om(N,u)$ cannot be compensated by regularity of $N$, and so
there seems no way to close the estimates as above for $(s,l)=(1,0)$.
Instead, we invoke the conservation laws with the weak compactness argument. This type of argument usually yields a weak result, typically without uniqueness. We can however obtain the strong well-posedness for small data as in Theorem \ref{thm1-1}, thanks to that both in the larger space $(s,0)$ with $s<1$ and in the smaller space $(1,l)$ with $l>0$.

Assume $(u_0,N_0)\in H^1\times L^2$. By Proposition \ref{small data scatter} there is $\e_0:=\e_1(1/2,0) \ll 1$ such that
if $\norm{(u_0,N_0)}_{H^{1/2}\times L^2} \leq \e_0$
then there is a unique global solution $(u,N)$ in $X^{1/2}\times Y^0$, satisfying
\begin{align}\label{eq:unispace}
  \norm{(u,N)}_{X^{1/2}\times Y^0}\leq C\e_0\ll 1.
\end{align}
Proposition \ref{regup Xs} implies that $(u,N)\in X^s\times Y^0$ for all $s\in[1/2,1)$.

Fix a sequence $\{(u_{0,n}, N_{0,n})\}\subset \cS(\R^4)$ such that $(u_{0,n},N_{0,n})\to (u_0,N_0)$ in $H^1\times L^2$ and $\norm{(u_{0,n},N_{0,n})}_{H^{1/2}\times L^2}\leq \e_0$.
By Proposition \ref{small data scatter}, for each $n$, there is a unique global solution $(u_n,N_n)$ satisfying \eqref{eq:unispace} and for all $1/2\le s<1$,
\begin{align}\label{eq:uniformbd}
  \sup_{n}\norm{(u_n,N_n)}_{X^{s}\times Y^0}<\infty.
\end{align}
Now we claim a uniform bound at the energy level:
\begin{equation} \label{Enb}
  \sup_{n,t}\|(u_n(t),N_n(t)) \|_{H^1\times L^2} < \I.
\end{equation}
By Proposition \ref{regup Xs}, we have $(u_n,N_n)\in X^{8}\times Y^{9}$ for all $n$, by which we can justify the conservation law $E_Z(u_n(t),N_n(t))=E_Z(u_{0,n},N_{0,n})$.
Using \eqref{eq:unispace} for $N_n$ together with the Sobolev inequality $\|u\|_{L^4_x}\lec \|\na u\|_{L^2_x}$ yields
\EQ{
 E_Z(u_n,N_n)=(1-O(\e_0))\norm{\nabla u_n}_2^2 + \|N_n\|_2^2/2,}
which, combined with the lower regularity bound \eqref{eq:uniformbd}, implies \eqref{Enb}.

Next we prove convergence $u_n(t)\to u(t)$ in $H^1_x$ as $n\to \infty$, locally uniformly in $\R$.
Take any convergent sequence $t_n\to t_\I$.
From Propositions \ref{small data scatter} and \ref{regup Xs}, we know that $u_n(t_n)\to u(t_\I)$ in $H^s_x$ for $s<1$, and $N_n(t_n)\to N(t_\I)$ in $L^2_x$.
From \eqref{Enb}, we have $\{u_n(t_n)\}_n$ is bounded in $H^1_x\subset L^4_x$, thus we get $u(t_\I)\in H^1$, $u_n(t_n)\to u(t_\I)$ weakly in $H^1_x$, and $|u_n(t_n)|^2\to |u(t_\I)|^2$ weakly in $L^2_x$.
Since $N_n(t_n)\to N(t_\I)$ strongly in $L^2_x$, we have
$\int N_n(t_n)|u_n(t_n)|^2dx \to \int N(t_\I)|u(t_\I)|^2 dx$, and so,
\EQ{ \label{weak conv EZ}
 E_Z(u(t_\I),N(t_\I)) \pt\leq \liminf_{n\to \infty} E_Z(u_n(t_n),N_n(t_n))
 \pr=\liminf_{n\to \infty} E_Z(u_{0,n},N_{0,n})=E_Z(u_0,N_0).}
By the time reversibility we get $E_Z(u(t_\I),N(t_\I))=E_Z(u_0,N_0)$.
Indeed, if there is a $t_0\in\R$ such that $E_Z(u(t_0),N(t_0))<E_Z(u_0,N_0)$, then we solve the Zakharov system with initial data $(u(t_0),N(t_0))$ at $t=t_0$.
By the uniqueness we get a contradiction.
Then the equality in \eqref{weak conv EZ} implies $\norm{\nabla u_n(t_n)}_{L^2}\to\norm{\nabla u(t_\I)}_{L^2}$, from which we conclude that $u_n(t_n)\to u(t_\I)$ strongly in $H^1_x$, and so the locally uniform convergence $u_n\to u$ in $C(\R;H^1_x)$.
Thus we obtain the unique global solution $(u,N)\in (C\cap L^\I)(\R;H^1\times L^2)$.
Note that the smoothness of the approximate solutions $(u_n,N_n)$ was used only to ensure the unique existence and the conservation law.
Now that we have them for the solutions in the energy space, we can apply the above argument to a sequence of initial data in $H^1\times L^2$, which implies continuous dependence of the initial data, locally uniformly in time.

By Propositions \ref{small data scatter} and \ref{regup Xs}, $(u,N)$ scatters to some $(u^+,N^+)$ in $H^s \times L^2$ for all $s<1$.
Since $u(t)\in L^\I(\R;H^1_x)$, we have $S(-t)u(t)\to u^+$ weakly in $H^1$ as $t\to +\I$. Since $|u(t)|^2$ is bounded in $(H^1_x)^2\subset B^1_{4/3}$, while $N(t)$ is vanishing in $B^{-1}_4$ as $t\to\I$ due to the scattering in $L^2_x\subset B^{-1}_4$, we have
\begin{align}\label{eq:decay1}
  \int N(t)|u(t)|^2dx \to 0 \quad (t\to +\I),
\end{align}
and so
\begin{align*}
  \|\na u^+\|_2^2 + \|N^+\|_2^2/2  \leq
  &\liminf_{t\to +\I} \|\na S(-t)u(t)\|_2^2 + \|W_\al(-t)N(t)\|_2^2/2\\
  =&\liminf_{t\to +\I}\|\na u(t)|_2^2 + \|N(t)\|_2^2/2\\
  =&\liminf_{t\to +\I}E_Z(u(t),N(t))=E_Z(u_0,N_0).
\end{align*}
To prove the equality above, we consider the final state problem.
Following the argument in Step 1, we fix a sequence $\{(u^+_n,N^+_n)\}\subset\cS(\R^4)$ such that $(u^+_n,N^+_n)\to(u^+,N^+)$ in $H^1\times L^2$.
Then by Proposition \ref{small data scatter}, we have a sequence of solutions $(\wt u_n,\wt N_n)\in X^{1/2}\times Y^0$ scattering to $(u^+_n,N^+_n)$ as $t\to\I$, which converges to $(u,N)$ in $X^{1/2}\times Y^0$ as $n\to\I$.
The regularity is upgraded to $X^s\times Y^l$ for all $(s,l)$ in Proposition \ref{regup Xs}.
As in Step 1, we have $\sup_{n,t}\norm{(\wt u_n,\wt N_n)}_{H^1\times L^2}<\I$, hence $u_n(t) \to u(t)$ weakly in $H^1_x$. Thus by \eqref{eq:decay1}
\begin{align*}
  E_Z(u(t),N(t))\leq \liminf_{n\to \I}E_Z(\wt u_n(t),\wt N_n(t))
  &=\liminf_{n\to \I} \|\na u_n^+|_2^2 + \|N_n^+\|_2^2/2 \\
  &= \|\na u^+|_2^2 + \|N^+\|_2^2/2.
\end{align*}
Hence we get
\[\lim_{t\to +\I}E_Z(u(t),N(t))= \|\na u^+\|_2^2 + \|N^+\|_2^2/2 \]
and so, $S(-t)u(t)\to u^+$ strongly in $H^1_x$, namely the scattering in $H^1_x$.

To show the continuity of the solution map in $L^\I_t(\R;H^1_x)$, it remains to prove $u_n(t_n)-u(t_n)\to 0$ in $H^1_x$ in the case $t_n\to\I$, for a sequence of solutions $(u_n,N_n)$ in the energy space such that $(u_n(0),N_n(0))\to (u(0),N(0))$ in $H^1\times L^2$.
Since $S(-t)u(t)\to u^+$ in $H^1_x$, it is equivalent to showing $S(-t_n)u_n(t_n)\to u^+$ in $H^1_x$.
We already know the $H^s_x$ convergence for $s<1$ as well as the weak convergence in $H^1_x$.
Then the strong convergence is equivalent to $\|u_n(t_n)\|_{H^1_x}\to\|u^+\|_{H^1}$. Since
\EQ{
 \|N_n(t_n)\|_{B^{-1}_4} \le \|N_n-N\|_{L^\I_tL^2_x}+\|N(t_n)\|_{B^{-1}_4} \to 0,}
we have $\int N_n(t_n)|u_n(t_n)|^2dx\to 0$, and so, as $n\to\I$,
\EQ{
 \pt \|\na u_n(t_n)\|_2^2 +\|N_n(t_n)\|_2^2/2 =E_Z(u_n(t_n),N_n(t_n))+o(1)
 \pr=E_Z(u_n(0),N_n(0))+o(1)=E_Z(u(0),N(0))+o(1)
 \pr=\|\na u^+\|_2^2 + \|N^+\|_2^2/2 +o(1). }
Since $\|u_n(t_n)\|_2\to \|u^+\|_2$ and $\|N_n(t_n)\|_2\to\|N^+\|_2$, the above implies the strong convergence of $S(-t_n)u_n(t_n)$ in $H^1_x$, and thus $u_n\to u$ in $L^\I_t(\R;H^1_x)$.
This completes the proof of Theorem \ref{thm1-1} in the case $(s,l)=(1,0)$.

\section{Ill-posedness at $(s,l)=(2,3)$}
In this section, we prove Theorem \ref{thm-ill}. The main point is that the multilinear estimates fail only for the boundary quadratic term coming from the initial data. Exploiting the dispersive smoothing, we can prove that the other terms are more regular if the initial data is localized in space.

\begin{proof}[Proof of Theorem \ref{thm-ill}]
First of all, for any initial data $(u_0,N_0)\in H^2\times H^3$, we have a unique local solution for $(s,l)$ in \eqref{sl range} satisfying $s\le 2$ and $l\le 3$, say $(u,N)\in(X^2\times Y^2)([0,T])$, by Propositions \ref{LWP large} and \ref{regup Xs}. In the Duhamel formula \eqref{eq:intN}, the first term on the right is obviously in $C(\R;H^3)$.
The integral terms are regular thanks to the high regularity. Indeed,
\EQ{
 \pt \|D|u|^2_{HH+\al L+L\al}\|_{L^1_tH^{3}_x} \lec \|u\|_{L^2_tB^2_4}^2,
 \pr \|D\ti\Om(Nu,u)\|_{L^1_tH^3_x} \lec \|N\|_{L^\I_tH^2_x}\|u\|_{L^2_tB^2_4}^2,}
and the same for $D\ti\Om(u,Nu)$.
To bound $D\ti\Om(u,u)$ in $H^3_x$, we use local smoothing for $u$, assuming that
\EQ{
 u_0\in W^{2,1}(\R^4)=\{f: \p^\al f\in L^1(\R^4) \mbox{ for } |\al|\leq 2\}.}
Then $S(t)u_0 \in C((0,\I);B^2_p)$ for all $p>2$ by the dispersive $L^p_x$ decay estimate for $S(t)$. Moreover, in the Duhamel formula \eqref{eq:intu} of $u$,
the terms except for $(Nu)_{LH\cdots}$ easily gain better regularity by
\EQ{
 \pt \|\Om(N,u)\|_{H^3_x} \lec \|N\|_{H^2_x}\|u\|_{H^2_x},
 \pr \|\Om(D|u|^2,u)\|_{L^1_tH^3_x} \lec \|u\|_{L^\I_tH^2_x}\|u\|_{L^2_tB^2_4}^2,
 \pr \|\Om(N,Nu)\|_{L^2_tB^4_{4/3}} \lec \|N\|_{L^\I_tH^2_x}^2\|u\|_{L^2_tB^2_4}.}
The remaining term is bounded in $C([0,T];B^2_3)$ by
\EQ{
 \norm{\int_0^tS(t-s)(Nu)_{LH\cdots}ds}_{B^2_3}
 \pt\lec \int_0^t|t-s|^{2/3}\|(Nu)_{LH\cdots}\|_{B^2_{3/2}}ds
 \pr\lec \int_0^t|t-s|^{2/3}\|N(s)\|_{L^6_x}\|u(s)\|_{H^2_x}ds.}
Gathering the above estimates, we obtain
$u \in C((0,T];H^2\cap B^2_3)$.
Since $B^2_3\subset L^\I$,
\EQ{
 \|D\ti\Om(u,u)\|_{H^3} \lec \|u\|_{H^2}\|u\|_{B^2_3},}
and plugging this into the above estimates for $N$, we deduce that
\EQ{
 N-W_\al(t)D\ti\Om(u_0,u_0) \in C((0,T];H^3_x),}
if $u_0\in H^2\cap W^{2,1}(\R^4)$.
Hence it suffices to find such a $u_0$ that $D\ti\Om(u_0,u_0)\not\in H^3_x$.
It is constructed in the next Lemma \ref{bad tiOm}. Then $N(t)\not\in H^3_x$ for all $0<t<T$, namely the instant exit or the latter part of the theorem.

Thanks to the high regularity, it is easy to translate it to non-existence. Indeed, if $(u,N)\in L^2((0,T);H^1\times H^3)$ then from the equation without the normal form,
\EQ{
 Nu\in L^1_tH^1_x \implies u\in C_tH^1_x\cap L^2_tB^1_4 \implies D|u|^2\in L^1_tL^2_x \implies N\in C_tL^2_x.}
In particular, $(u,N)$ belongs to the uniqueness class at $(s,l)=(1/2,0)$.
Hence it should be identical with the exiting solution obtained above, satisfying $N(t)\not\in H^3$ for all $t\not=0$, contradicting $N\in L^2_t((0,T);H^3_x)$.
\end{proof}
It remains to prove the failure of the bilinear estimate:
\begin{lem} \label{bad tiOm}
There is a radial $u\in (H^2\cap W^{2,1})(\R^4)$ satisfying $D\ti\Om(u,u)\not\in H^3(\R^4)$.
\end{lem}
This failure of bilinear estimate comes from that $H^2(\R^4)$ is not an algebra, but we should be careful about cancellation in the nonlinearity. In fact, the proof below implies that $D\ti\Om(u,u)$ is bounded in $H^3$ for real-valued or purely imaginary $u\in H^2$.

\begin{proof}
Modulo a bounded operator, the symbol of $D\ti\Om$ can be approximated
\EQ{
 \frac{\al|\x|}{|\x-\y|^2-|\y|^2\mp\al|\x|}=\frac{\al}{|\x|}+\frac{\al(2\x\cdot\y\pm\al|\x|)}{|\x|(|\x-\y|^2-|\y|^2\mp\al|\x|)}}
in the $XL$ frequency, while in the $LX$ frequency,
\EQ{
 \frac{\al|\x|}{|\x-\y|^2-|\y|^2\mp\al|\x|}=\frac{-\al}{|\x|}+\frac{\al(2\x\cdot(\x-\y)\pm\al|\x|)}{|\x|(|\x-\y|^2-|\y|^2\mp\al|\x|)},}
where the second terms are $O(|\x|^{-2}\LR{Low})$ for $|\x|\gg 1$, and so bounded $H^2\times H^2\to H^3$ for high frequency.
Hence it suffices to construct $u \in H^2\cap W^{1,2}$ such that $\supp\hat u(\x)=0$ for $|\x|\lec 1$ and
\EQ{ \label{tiOm main}
 (u \bar u)_{HL} - (u \bar u)_{LH} \not \in H^2(\R^4).}
Indeed, this is necessary and sufficient for $D\ti\Om(u,u)\not\in H^3$ under the condition of $\supp\hat u$. Note that the left side is simply zero if $u(\R^4)\subset\R$ or $iu(\R^4)\subset\R$. The remaining is the anti-symmetric part, which can be expanded by putting $u=v+iw$
\EQ{
 (u \bar u)_{HL} - (u \bar u)_{LH}=2i[(wv)_{HL}-(wv)_{LH}].}
Now it is easy to avoid the cancellation considering the form
\EQ{
 v=\sum_{j>J}a_j \fy_j, \pq w=\sum_{j> J}b_j \fy_j,
 \pq \fy_j(x)=\fy(2^jx),}
where $J\gg\log\al$, $\{a\},\{b\}\subset[0,\I)$, and $\fy\in\cS(\R^4;\R)$ is a non-zero real-valued radial function satisfying
\EQ{
 0\le \hat\fy \le 1, \pq \supp\hat\fy\subset\{||\x|-1|\ll 1\}.}
Put $c:=\fy(0)>0$. Injecting the above ansatz expands the bilinear form
\EQ{ \label{vw in tiOm}
 (vw)_{HL}-(vw)_{LH}=\sum_{j>J}\sum_{k>J}^{j-K}(a_jb_k-a_kb_j)\fy_j\fy_k.}
Since $\F(\fy_j\fy_k)$ is supported around $|\x|=2^j$,
\EQ{
 \|\eqref{vw in tiOm}\|_{H^2}^2 \sim \sum_{j> J}\|2^{2j}\sum_{k\le j-K}(a_jb_k-a_kb_j)\fy_j\fy_k\|_2^2.}
Imposing a support condition on $\{a\},\{b\}$
\EQ{
 \supp a \cap \supp b = \emptyset,}
we can decouple the above as
\EQ{
 \|\eqref{vw in tiOm}\|_{H^2}^2 \sim \sum_{j> J}\|2^{2j}\sum_{k\le j-K}a_jb_k\fy_j\fy_k\|_2^2 + \sum_{j> J}\|2^{2j}\sum_{k\le j-K}b_ja_k\fy_j\fy_k\|_2^2.}
By rescaling $x\mapsto 2^{-j}x$, and using $\fy(2^{k-j}x)=c+O(|2^{k-j}x|)$, the $L^2_x$ norm is approximated by
\EQ{
 \pt\|2^{2j}\sum_{k\le j-K}a_jb_k\fy_j\fy_k\|_{L^2_x}
 = \|a_j\fy(x)\sum_{k\le j-K}b_k\fy(2^{k-j}x)\|_{L^2_x}
 \pr\ge c|a_j|\|\fy\|_{L^2_x}\sum_{k\le j-K}b_k-C\|a_jx\fy(x)\|_{L^2_x}\sum_{k\le j-K}b_k2^{k-j}.}
Fix $\te\in(1/2,3/4)$ and let
\EQ{
 a_j=\CAS{ j^{-\te} &J< j\text{ is even}\\ 0 &\text{otherwise},}
 \pq b_j=\CAS{ j^{-\te} &J< j\text{ is odd}\\ 0 &\text{otherwise}.}}
Then for $j>K+J$,
\EQ{
 \sum_{k\le j-K}b_k \sim (j-K)^{1-\te}, \pq \sum_{k\le j-K}b_k2^{k-j} \lec 2^{-K},}
and so
\EQ{
 \|\eqref{vw in tiOm}\|_{H^2} \gec \|j^{-\te}(j-K)^{1-\te}\|_{\ell^2(j>J+K)}-C2^{-K}\|j^{-\te}\|_{\ell^2(j>J)}=\I,}
since $-\te<-1/2<1-2\te$. Also we have
\EQ{
 \|u\|_{H^2}\lec\|j^{-\te}\|_{\ell^2(j>J)}<\I,
 \pq \|u\|_{W^{2,1}} \lec \|2^{-2j}j^{-\te}\|_{\ell^1(j>J)}<\I.}
Thus we have obtained a desired example $u\in H^2\cap W^{2,1}$.
\end{proof}

\subsection*{Acknowledgment}
Z.\ Guo is supported in part by NNSF of China (No.\ 11371037), Beijing
Higher Education Young Elite Teacher Project (No.\ YETP0002), and Fok
Ying Tong education foundation (No.\ 141003).

S.\ Herr was supported by the German Research Foundation, CRC 701.

Part of this research was carried out while the authors participated in the program ``Harmonic Analysis and Partial Differential Equations'' at the Hausdorff Research Institute for
Mathematics in Bonn.

\bibliographystyle{amsplain}

\end{document}